\documentclass[review]{elsarticle}

\usepackage{lineno,hyperref}
\usepackage{moreverb}
\usepackage{graphicx, wrapfig}
\usepackage{algorithm}
\usepackage{algpseudocode}
\usepackage{color}
\usepackage{amsmath}
\usepackage{txfonts}
\usepackage{afterpage}
\usepackage[margin=1in]{geometry}
\usepackage{mathrsfs}
\usepackage{tikz}
\usepackage{pstricks}
\usepackage{diagbox}
\usepackage{arydshln}
\usepackage{dsfont}
\usepackage{smartdiagram}
\usepackage{todonotes}
\usepackage{verbatim}
\usepackage{smartdiagram}
\usepackage{bm}
\usepackage{commath}
\usepackage{bbm}
\usepackage{mathrsfs}
\usepackage{dsfont}
\usepackage{extarrows}
\usepackage{subfigure}
\usepackage{multirow}
\usepackage{graphicx}
\usepackage{array}

\DeclareMathOperator*{\argmin}{arg\,min}
\DeclareMathOperator*{\argmax}{arg\,max}

\newcommand{\bx}{x}

\newcommand{\bL}{\bm{L}}
\newcommand{\bxi}{\xi}

\newcommand{\calL}{\mathcal{L}}
\newcommand{\calD}{\mathcal{D}}

\newcommand{\bdobs}{\bm{d}_{\text{obs}}}
\newcommand{\dobs}{d_{\text{obs}}}
\newcommand{\bxobs}{\bm{x}_{\text{obs}}}
\newcommand{\bxobsij}{\bm{x}_{\text{obs}}^{(i,j)}}
\newcommand{\pC}{\mathfrak{c}}
\newcommand{\frakN}{\mathfrak{N}}
\newcommand{\NC}{N}
\newcommand{\calB}{\mathcal{B}}
\newcommand{\bbR}{\mathbb{R}}

\newcommand{\checkE}{\check{\mathbb{E}}}
\newcommand{\tauij}{\tau^{(i,j)}}
\newcommand{\hij}{h^{(i,j)}}
\newcommand{\gij}{g^{(i,j)}}
\newcommand{\tilgij}{\widetilde{g}^{(i,j)}}
\newcommand{\tila}{\widetilde{a}}
\newcommand{\tilpsi}{\widetilde{\psi}}
\newcommand{\calDi}{\mathcal{D}^{(i)}}
\newcommand{\calDj}{\mathcal{D}^{(j)}}
\newcommand{\hatxi}{\widehat{\xi}}
\newcommand{\hata}{\widehat{a}}
\newcommand{\brea}{\breve{a}}
\newcommand{\atruth}{a_{\text{truth}}}

\newcommand{\obsvar}{\sigma_{\text{obs}}}
\newcommand{\tolKL}{\delta_{\text{KL}}}
\newcommand{\nij}{n^{(i,j)}}
\newcommand{\Lambdaij}{\Lambda^{(i,j)}}

\newcommand{\Npara}{N_{\xi}}

\newcommand{\checkV}{\check{\mathbb{V}}}
\newcommand{\ccalD}{\overline{\calD}}
\newcommand{\ccalDi}{\overline{\calDi}}

\def\covl{L}
\newcommand{\reint}{\epsilon_{\text{int}}}
\newcommand{\restate}{\epsilon_{\text{state}}} 

\usetikzlibrary{calc,trees,positioning,arrows,chains,shapes.geometric,%
    decorations.pathreplacing,decorations.pathmorphing,shapes,%
    matrix,shapes.symbols}

\tikzset{
>=stealth',
  punktchain/.style={
    rectangle, 
    rounded corners, 
    draw=black, very thick,
    text width=10em, 
    minimum height=2em, 
    text centered, 
    on chain},
  line/.style={draw, thick, <-},
  element/.style={
    tape,
    top color=white,
    bottom color=blue!50!black!60!,
    minimum width=6em,
    draw=blue!40!black!90, very thick,
    text width=10em, 
    minimum height=3.5em, 
    text centered, 
    on chain},
  every join/.style={->, thick,shorten >=1pt},
  decoration={brace},
  tuborg/.style={decorate},
  tubnode/.style={midway, right=2pt},
}

\tikzstyle{arrow} = [thick,->,>=stealth]
\newtheorem{thm}{Theorem}
\newtheorem{Def}{Definition}

\newdefinition{rmk}{Remark}
\newproof{pf}{Proof}

\newcommand{\comm}[1]{}

\makeatletter
\def\ps@pprintTitle{%
   \let\@oddhead\@empty
   \let\@evenhead\@empty
   \let\@oddfoot\@empty
   \let\@evenfoot\@oddfoot
}
\makeatother

\modulolinenumbers[5]

\journal{Journal of Computational Physics}
\begin{document}

\begin{frontmatter}
		
		\title{Domain-decomposed Bayesian inversion based on local Karhunen-Lo\`{e}ve expansions}
		
		%%% Group authors per affiliation:
		
		\author[ShanghaiTechAddress]{Zhihang Xu}
		\ead{xuzhh@shanghaitech.edu.cn}
		
    \author[ShanghaiTechAddress]{Qifeng Liao\corref{mycorrespondingauthor}}
		\cortext[mycorrespondingauthor]{Corresponding author}
		\ead{liaoqf@shanghaitech.edu.cn}

		\author[JinglaiLiAddress]{Jinglai Li}
		\ead{j.li.10@bham.ac.uk}
		
		\address[ShanghaiTechAddress]{School of Information Science and Technology, ShanghaiTech University, Shanghai 201210, China}
		\address[JinglaiLiAddress]{The School of Mathematics, University of Birmingham, Birmingham B15 2TT, UK}

		\begin{abstract}    
    In many Bayesian inverse problems the goal is to recover a spatially varying random field.
    Such problems are often computationally challenging especially when the forward model is governed by complex partial differential equations (PDEs).  
    The challenge is particularly severe when the spatial domain is large and the unknown random field needs to be represented by a high-dimensional parameter.
    In this paper, we present a domain-decomposed method to attack the dimensionality issue and the method decomposes the spatial domain and the parameter domain simultaneously.
    On each subdomain, a local Karhunen-Lo\`{e}ve (KL) expansion is constructed, and a local inversion problem is solved independently in a parallel manner, and more importantly, in a lower-dimensional space. After local posterior samples are generated through conducting Markov chain Monte Carlo (MCMC) simulations on subdomains, a novel projection procedure is developed to effectively reconstruct the global field. In addition, the domain decomposition interface conditions are dealt with an adaptive Gaussian process-based fitting strategy. Numerical examples are provided to demonstrate the performance of the proposed method.
		\end{abstract}
		
		
	\end{frontmatter}

\section{Introduction}
%% Definition/Main task of inverse problems
%% Its applications in various fields 
%% Its difficulties: complex system (that is why we use DD), high- dimensionality and ill-poseness.
Many real world inverse problems involve forward models governed by partial differential equations (PDEs), and
in these problems often the primary task is to recover spatially varying unknown parameters from noisy and incomplete observations.
Such problems are ubiquitous in various scientific areas, including geosciences \cite{springer2021efficient}, climate prediction \cite{solonen2012efficient},
seismic inversion \cite{martin2012stochastic} and remote sensing \cite{haario2004markov}.
The Bayesian inference \cite{stuart2010inverse,
wang04bayesian,efendiev2006preconditioning,marzouk2007stochastic,lieberman2010parameter,elmoselhy12bayesian,tan13computational,li2014adaptive, cui16dimension,zahm2022certified}  
has become an important tool for solving such problems, largely due to its ability to quantify the uncertainty in the solutions obtained. 

While the Bayesian methods  are conceptually straightforward,  their applications to the 
aforementioned PDE-involved inverse problems can be extremely challenging, where 
a major difficulty lies in the computational aspect. 
As is well known, in most practical problems, the posterior distributions are analytically intractable, and are often 
computed with sampling methods.
One of the most popular methods in this context is the Markov chain Monte Carlo (MCMC) simulation \cite{robert2013monte}. 
The major limitation associated with MCMC as well as other sampling methods is that they typically require a very large number of 
evaluations of the forward model, which can be prohibitively costly for our problems, as the PDE-involved forward model is computationally intensive. 
While considerable efforts have been devoted to reducing the computational cost, e.g., \cite{li2014adaptive,cui2015data,chen15sparse, jiang17multiscale,liao2019adaptive}, many challenges remain in
applying the Bayesian methods for the PDE-involved inverse problems.  
Among them, the dimensionality issue is one of the most frequently encountered difficulties in these problems. 

To conduct Bayesian inference, one first needs to parametrize the spatially varying unknown parameter (in Bayesian inference it is typically modelled 
as a random field) as a finite-dimensional parameter. 
%Parameterization of the spcially varying random field can represent an infinite field through a set of finite random variables and provide a way to characterizes %the parameter in a lower-dimensional space.
Existing methods for doing so include the Karhunen-Lo\`{e}ve (KL) expansion 
\cite{ghanem2003stochastic,li2015note}, 
wavelet-based parameterization \cite{ellam2016bayesian}, 
and parameterization techniques based on deep generative models (DGM) \cite{xia2021bayesian}. 
In this paper, we focus on the KL expansion since it is optimal in the mean squared error sense with respect to the number of random variables in the representation. 
In many practical problems, especially those with large spatial domains,  often a large number of KL modes
are needed to represent the unknown field, leading to a very high-dimensional inference problem.
 The  primary goal of this work is to address this issue and 
 reduce the dimensionality of the inverse problems using a domain decomposition (DD) approach. 
In particular, we perform domain decomposition  over the spatial domain and the parameter space simultaneously. 
The resulting method enables parallelization and thus facilitates efficient sampling in a much lower dimensional parameter space.

In general, domain decomposition for uncertainty quantification and inverse problems gains a lot of interests, and related methods are actively developed. In \cite{chen2015local}, local polynomial chaos expansions based on domain decomposition are proposed for solving PDEs with high dimensional random inputs. In \cite{liao2015domain}, we provide a domain-decomposed uncertainty quantification approach based on importance sampling. Efficient methods to compute  dominant KL terms through domain decomposition and the corresponding accelerated Monte Carlo sampling procedures are presented in \cite{contreras2018parallela,contreras2018parallelb}. Domain decomposition methods for solving nonlinear transient inverse heat conduction problems are studied in \cite{khajehpour13domain}. In \cite{jagtap2020conservative, jagtap2020extended,shukla2021parallel}, domain decomposition methods with physics-informed neural networks are addressed for forward and inverse problems.

In this work, we focus on domain decomposition for Bayesian inversion, and the main contributions of this work are as follows. The first is effective local representation for priors. It is known that when the unknown fields have complex structures, the corresponding global priors need to have short correlation lengths to give effective inference results, which requires high-dimensional global parameterization. In the method proposed in this work, relative correlation lengths are increased along with decomposing a global spatial domain into small subdomains, such that low-dimensional parameters can approximate complex priors well. The second is efficient forward model evaluation procedures.  As discussed above, main computational costs of sampling based inference methods are caused by repetitively evaluating expensive forward models, especially for models governed by PDEs.  High-fidelity numerical schemes can give accurate predictions for the outputs of these PDEs, e.g., the finite element methods with a posteriori error bounds \cite{AINSWORTH19971,elman14finite}, but they can be expensive, as they require a large number of degrees of freedom when the underlying model is complex. As the global spatial domain is decomposed in our setting, the finite element degrees of freedom on local subdomains are significantly smaller than those for the global domain, and therefore evaluating each local model is clearly cheaper than evaluating the global model. The third is a new reconstruction approach for the global unknown field. Once local inversions are conducted, directly stitching local fields to approximate global unknown fields can give misleading information on domain decomposition interfaces. For this purpose, our new approach conducts projection of the local inference results over the space spanned by the global KL modes, which gives an effective approximation for the global true field. Lastly, to give proper interface conditions, Gaussian process (GP) models for interface treatments are built with an active learning procedure.

The rest of this paper is organized as follows. 
Section \ref{sec:problem_setup} sets the problem, where the standard MCMC procedure and the KL expansion are reviewed. 
In Section \ref{sec:dd_para}, we discuss the KL expansion on local subdomains, and give our reconstruction procedure for the global input fields. In Section \ref{sec:dd_mcmc},  our Gaussian process interface treatments are discussed, and our overall domain-decomposed Markov chain Monte Carlo (DD-MCMC) algorithm is presented. Numerical results are discussed in Section \ref{sec:numerical}. 
Section \ref{sec:conclusion} concludes the paper.

%In section \ref{sec:problem_setup}, we introduce our problem setup and  discuss the smooth local parameterization is  in section \ref{sec:dd_para}. Section \ref{sec:dd_mcmc} presents the domain-decomposed Markov chain Monte Carlo (DD-MCMC) algorithm in detail, and section \ref{sec:numerical} provides numerical examples to demonstrate the performance of the proposed method. Finally, concluding remarks are offered in section \ref{sec:conclusion}. 

\section{Problem setup}
\label{sec:problem_setup}
In this section, we briefly review the general description of Bayesian formulation for inference and detailed settings for KL expansion parameterization for PDEs with random inputs.

\subsection{Bayesian inverse problem}
Letting $\xi$ denote a $\Npara$-dimensional parameter of interest and $\bdobs\in \bbR^{n} (n\ll \Npara)$ denote  $n$-dimensional observed data, we want to estimate $\xi$ from $\bdobs$.
We assume that there exists a forward model $F$ that maps the unknown parameter $\xi$ to the data $\bdobs$:
\begin{equation}
\label{eq:forward_model}
    \bdobs = F(\xi)+ \epsilon\,,
\end{equation}
where $\epsilon\in\bbR^{n}$ denotes the random observation noise and 
its probability density function is denoted by $\pi_{\epsilon}(\epsilon)$. 
Then, the likelihood function which characterizes the relationship between observations and the forward model can be defined as 
\begin{eqnarray}
\bL(\bdobs|\xi) = \pi_{\epsilon}(\bdobs - F(\xi))\,.
\label{eq_likelihood}
\end{eqnarray}
In this paper, the noise $\epsilon$ is assumed to be Gaussian with zero mean and a diagonal covariance matrix  
$\obsvar^2 \bm{I}_{n}$, where $\obsvar >0$ is the standard deviation and $\bm{I}_n$ is the identity matrix with size $n\times n$.
The likelihood function is then proportional to the data-misfit functional
$\eta(\xi;\bdobs) := \frac{1}{2\sigma^2}\|\bdobs - F(\xi)\|_2^2$ where $\| \cdot\|_2$ denotes the Euclidean norm, i.e.,
\begin{equation*}
	\label{eq:likelihood}
	\bL(\bdobs|\xi) \propto \exp(-\eta(\xi;\bdobs)).
\end{equation*}
Given a prior distribution $\pi_0(\xi)$ of $\xi$ which reflects the knowledge of the parameter before any measurements, based on the Bayes' rule,  
the posterior distribution of $\xi$ can be written as
\begin{equation}
	\label{eq:bayes}
	\pi(\xi|\bdobs) = \frac{ \overbrace{\bL(\bdobs|\xi)}^{\text{likelihood}}
		\overbrace{\pi_0(\xi)}^{\text{prior}}
	}
	{
		\underbrace{
			\pi(\bdobs)}_{\text{evidence}}}
	\propto \bL(\bdobs|\xi) \pi_0(\xi)
	\,,
\end{equation}
where the evidence $\pi(\bdobs)$ in \eqref{eq:bayes} is usually viewed as a normalization constant for a well-defined probability distribution. 
The posterior distribution is usually analytically intractable, and therefore sampling methods including Markov chain Monte Carlo methods (MCMC) \cite{robert2013monte} are widely used.
A Markov chain is a sequence of samples where the next state only depends on the previous state, which is known as the Markov property, and the move from the current state towards the next state is defined through some transition operator.
The MCMC method constructs a Markov chain, of which the equilibrium distribution (also known as the stationary distribution) is set to the target distribution. 
In the context of Bayesian inversion, the target distribution is the posterior distribution. 
To ensure the convergence towards the target distribution, the detailed balance condition has to be satisfied.
To generate samples of the posterior distribution, we consider the standard Metropolis-Hastings (MH) \cite{metropolis1953equation, hastings1970monte} 
algorithm, which proceeds as follows. 
Starting from a randomly chosen initial state, for the $s$-th state $\xi^{s}$, a candidate state $\xi^{\star}$ is drawn from some proposal distribution $Q(\cdot|\xi^{s})$, 
and then the candidate state is accepted with the probability of an
acceptance rate denoted by $\alpha(\xi^{\star},\xi^{s})$. 
The proposal distribution and the acceptance probability define the transition operator, i.e., $h(\xi^{\star},\xi^{s}) = Q(\xi^{\star}|\xi^{s}) \alpha(\xi^{\star},\xi^{s}).$
The detailed balance condition is given through the transition operator,  
\[
\pi(\xi^s|\bdobs) h(\xi^{s},\xi^{\star}) = \pi(\xi^{\star}|\bdobs) h(\xi^{\star},\xi^{s})\,.
\]
To guarantee that the detailed balance condition is satisfied, the acceptance probability can be defined as 
\[
\alpha(\xi^{\star},\xi^{s}) = \min\left\{1,
\frac{Q(\xi^{s}|\xi^{\star})\bL(\bdobs|\xi^{\star})\pi_0(\xi^{\star})}{Q(\xi^{\star}|\xi^{s})\bL(\bdobs|\xi^{s})\pi_0(\xi^{s})}\right\}\,.
\]
Details of the MH approach is summarized in Algorithm \ref{alg:MHMCMC}, where $\NC$ is a given number of posterior samples to generate. 

\begin{algorithm}[!ht]
    \caption{The standard MH algorithm}
    \label{alg:MHMCMC}
    \begin{algorithmic}[1]
        \Require{Forward model $F(\xi)$, observation data $\bdobs$.}
        \State Generate an initial state $\xi^{1}$. 
    \For{$s=1,\ldots,\NC-1$}
    \State Draw $\xi^{\star}$ from a proposal distribution $Q(\cdot|\xi^{s})$. 
    \State 
    Compute the acceptance ratio 
    	\[
\alpha = \min\left\{1,
\frac{Q(\xi^{s}|\xi^{\star})\bL(\bdobs|\xi^{\star})\pi_0(\xi^{\star})}{Q(\xi^{\star}|\xi^{s})\bL(\bdobs|\xi^{s})\pi_0(\xi^{s})}\right\}\]
where $\pi_0$ is a given prior distribution and the likelihood  $\bL$
defined in \eqref{eq_likelihood} requires the forward model \eqref{eq:forward_model}.
    \State Draw $\rho$ from a uniform distribution $\rho \sim \mathcal{U}[0,1]$.
    \If {$\rho<\alpha$}
    \State
    Accept the proposal state, i.e., let $\xi^{s+1} = \xi^\star$.
    \Else
    \State Reject the proposal state, i.e., let $\xi^{s+1}= \xi^{s}$.
    \EndIf
    \EndFor
    \Ensure Posterior samples $\{\xi^{s}\}_{s=1}^{\NC}$.
    \end{algorithmic}
\end{algorithm}

\subsection{PDEs with random inputs and parameterization}
This section presents the detailed settings of the forward model considered in the paper.
Let $\mathcal{P}=(\Omega, \Sigma_{\Omega}, \mu_{\Omega})$ be a probability space, where $\Omega$ is the set of events, $\Sigma_{\Omega}$ is a sigma-algebra over $\Omega$ and $\mu_{\Omega}$ is a probability measure.
We denote the expectation operator for a function $\mathcal{F}(\cdot)$ as 
\[
\mathbb{E}[\mathcal{F}] = \int_{\Omega}\mathcal{F}(\omega) \dif \mu_{\Omega}(\omega)\,,
\]
and denote $L_2(\Omega)$ the space of second-order random variables, i.e., $L_2(\Omega):= \{\mathcal{F}|\mathbb{E}[\mathcal{F}^2]<+\infty \}$.
Moreover, $\calD \subset \mathbb{R}^{N_D}$ ($N_D=1,2,3$) denotes a physical domain which is  bounded, connected and with a polygonal boundary $\partial \calD$,
and $x\in \calD$ denotes a spatial variable. 
The space of square integrable functions is defined as  
$L_2(\calD) :=\{ \mathcal{F}| \int_{\calD}\mathcal{F}^2 < +\infty \}$, 
and the corresponding inner product is defined as $\langle \mathcal{F}(x),\mathcal{G}(x) \rangle_{\calD}:= \int_{\calD} \mathcal{F}(x)\mathcal{G}(x)\dif x$ for any $\mathcal{F}$ and $\mathcal{G}$ belonging to $L_2(\calD)$.
For any $\mathcal{F}\in L_2(\calD)$, the norm induced by the inner product is defined by 
\[
\|\mathcal{F}\|_{\calD}^2 = \langle \mathcal{F},\mathcal{F}\rangle_{\calD} =\int_{\calD} \mathcal{F}(x)^2\dif x.
\]
%{\color{red} %use math definiton
%Then, $L_2(\calD,\Omega)$ denotes the space of real-valued second-order processes $\kappa: \calD \times \Omega \to \mathbb{R}$, i.e., $L_2(\calD,\Omega):= \{\kappa|\kappa(x,\cdot)\in L_2(\calD)\text{ and } \kappa(\cdot, \omega) \in L_2(\Omega)\}$.
%}
The physics of problems considered 
are governed by 
a PDE over the spatial domain $\calD$ and
boundary conditions on the boundary $\partial \calD$,
which are stated as: 
find
$v(x,\omega)$ mapping $\calD\times \Omega$ to $\mathbb{R}$, such that
%\begin{equation}
\begin{subequations}    
\begin{align}
	&\calL(x,v; \kappa(x,\omega)) = f(x)\,,\quad x \in \calD\,,\\	
	&\calB(x,v; \kappa(x,\omega)) = h(x)\,,\quad x \in \partial \calD\,, 
\end{align}
\label{eq:original_pde}
\end{subequations} 
%\end{equation}
where $\calL$ is a differential operator and $\calB$ is a boundary condition operator, both of which are dependent on a random field $\kappa(x,\omega)$. 
Here $f$ is the source term 
and $h$ specifies the boundary condition. 

Generally, the random field $\kappa(x,\omega)$ is infinite-dimensional and needs to be parameterized.
As the truncated Karhunen-Lo\`{e}ve (KL) expansion is an optimal representation of random processes in the mean squared error sense, we focus on this expansion.
Letting $a_0(x)$ be the mean function of $\kappa(x,\omega)$, the covariance function $C(x,y):\calD\times\calD\to \bbR$ is defined as 
\[
	C(x,y) = \mathbb{E}[(\kappa(x,\omega) - a_0(x))(\kappa(y,\omega) - a_0(y))]\,,\quad x,y \in \calD\,.
\]
We can express the covariance function as $C(x, y) = \sigma(x) \sigma(y) \rho(x,y)$, where $\sigma: \calD \to \mathbb{R}$ is the standard deviation function of the random field and $\rho: \calD\times \calD \to [-1,1]$ is its autocorrelation coefficient function.
Let $\{\lambda_r,\psi_r(x)\}_{r=1}^{\infty}$ be the eigenvalues and the associated orthonormal eigenfunctions of the covariance function, that is, they satisfy 
\begin{equation}
\label{eq:global_eigenpair_1}
    \int_{\calD} C(x,y) \psi_r(x)\dif x = \lambda_r\psi(y)\,,\quad r=1,2,\ldots,
    \quad x,y\in \calD\,,
\end{equation}
and
\begin{equation}
\label{eq:global_eigenpair_2}
    \int_{\calD} \psi_r(x)\psi_t(x)\dif x = \delta_{rt}\,,
\end{equation}
where $\delta_{rt}$ denotes the Kronecker delta, and here we assume that the eigenvalues are ordered in decreasing magnitude.
By Mercer's Theorem, the covariance function has the following spectral decomposition, 
\[
C(x,y)  =\sum_{r=1}^{\infty} \lambda_r \psi_r(x) \psi_r(y)\,.
\]
According to the decomposition, it can be seen that 
$
 \sum_{r=1}^{\infty} \lambda_r = \int_{\calD} C(x,x)\dif x\,.
$
Based on the eigen-decomposition of the covariance function, the KL expansion provides a representation in terms of infinite number of random variables,
\begin{equation}
\label{eq:kle_inf}
	\kappa(x,\omega) =
	a_0(x) + \sum_{r=1}^{\infty}\sqrt{\lambda_r} \psi_r(x) \xi_r(\omega)\,,\quad x\in \calD\,,
\end{equation}
where $\{\xi_r(\omega)\}_{r=1}^{\infty}$ are uncorrelated random variables which control the randomness of the filed.
For a given random field  $\kappa(x,\omega)$, 
the corresponding random variables can be given via the orthonormality of eigenfunctions, 
\begin{equation*}
	%\label{eq:kle_xi}
	\xi_r(\omega) = \frac{1}{\sqrt{\lambda_r}} \int_{\calD} [\kappa(x,\omega)) - a_0(x)] \psi_r(x) \dif x\,,\quad r=1,2,\ldots\,,
\end{equation*}
and satisfy $\mathbb{E}[\xi_r] =0 $ and $\mathbb{E}[\xi_r \xi_t] =\delta_{rt}.$
For practical implementations, \eqref{eq:kle_inf} can be truncated with a finite number of terms such that the leading-order terms are maintained, 
\begin{equation}
\label{eq:kle_global}
\kappa(x,\omega) \approx a(x, \xi(\omega)) = a_0(x) + \sum_{r=1}^{d}\sqrt{\lambda_r} \psi_r(x) \xi_r(\omega) \,,\quad x\in \calD\,,
\end{equation}
where $\xi(\omega):=[\xi_1(\omega),\ldots,\xi_d(\omega)]^T$.
In this paper, we refer to $a(x,\xi(\omega)) - a_0(x)$ as the \textit{centralized} random field of $a(x,\xi(\omega))$.
The truncation level $d$ depends on the decay rate of eigenvalues which is related to the correlation length of the random field.
Usually, we select $d$ such that at least $\tolKL$ (a given threshold) of the total variance is captured, i.e., 
\begin{eqnarray}
\left(\sum_{r=1}^d \lambda_r\right)/(|\calD|\sigma^2)>\tolKL, \label{eq_tolKL}
\end{eqnarray}
where $|\calD|$ denotes the area of the domain $\calD$.
The prior distribution of $\xi$ is denoted by $\pi_0(\xi)$, of which the support is denoted by $I_{\xi}\subset \mathbb{R}^d$.
For a continuous covariance  function, the truncated KL expansion converges in the mean square sense uniformly \cite{le2010spectral} on $\calD$, i.e., 
\[
\lim_{d\to \infty} \sup_{x\in \calD} 
\mathbb{E}\left(
\kappa(x,\omega) - a_0(x) - \sum_{r=1}^d \sqrt{\lambda_r}\psi_r(x)\xi_r(\omega)
\right)^2  =0\,.
\]

After the above parameterization procedure over the random field, the original governing equation \eqref{eq:original_pde} is then transformed into the following finite-dimensional parameterized PDE system: find $u(x,\xi)$ mapping $\calD\times I_{\xi}$ to $\mathbb{R}$ such that
\begin{subequations}
\begin{align}
	&\calL(x,u; a(x,\xi)) = f(x)\,,\quad x \in \calD\,,\\	
    &\calB(x,u; a(x,\xi)) = h(x)\,,\quad x \in \partial \calD\,.
\end{align}
\label{eq:parameterized_pde}
\end{subequations}
Through specifying an observation operator $\pC$, e.g., 
taking solution values at several grid points, 
we write the overall forward model as
$F(a(x,\xi)):=\pC(u(x,\xi))$. 
%and consequently, the observational data is given by $\bdobs=F(\xi)+\epsilon$ where $\xi=a$.
% motivation for DD-inverse

As discussed in detail in \cite{chen2015local},  
for a given random field with correlation length $L_a$, 
the decay rate of the eigenvalues (see \eqref{eq:global_eigenpair_1}) depends on
the relative correlation length, i.e., $L_{a,\calD} := L_a /L_{\calD}$, where $L_{\calD}$ is the diameter of the physical domain $\calD$. 
It is shown that long correlation lengths lead to fast decay of eigenvalues, and vice versa.
So, when the  correlation length of the global field $\kappa(x,\omega)$ is small, its parameterization over the global domain $\calD$ 
can be high-dimensional (i.e., $d$ in \eqref{eq:kle_global} is large).
To result in a low-dimensional parameterization, we next decompose the physical domain into small subdomains, and the relative correlation length then becomes larger.

%As discussed previously, when the random field has a short correlation length, the parameterized input $a(x,\xi)$ still needs a large number of terms \cite{chen2015local, hou2017exploring}, and consequently, such inverse problem is in general high-dimensional and ill-posed. To explore the high-dimensional parameter space, MCMC may require tremendous times of repeated evaluations of the PDE system \eqref{eq:parameterized_pde} for different values of random inputs, which can be computationally demanding. In this work, decomposition over the spatial domain provides a viable route to alleviate the two major bottlenecks. After decomposition, the original problem is transformed into a sequence of sub-problems. Locally, each sub-problem can be considered independently and more importantly, in a lower-dimensional space. In this way, MCMC can achieve better sampling behavior in local problems, i.e., less mixing time and more effective samples given fixed computational resources. Moreover, by considering smaller computational region, the cost of the local forward evaluation is also drastically reduced. 

\section{Domain-decomposed Parameterization}
\label{sec:dd_para}
In this section, we first discuss settings for domain decomposed local problems with KL expansion parameterization posed on subdomains.  After that, based on realizations of local KL expansions, a new procedure to reconstruct global permeability fields is presented. These reconstructed global fields are called the assembled fields, and they are shown to be the projections of local fields to the space spanned by global eigenfunctions in KL expansion. 

\subsection{Local KL expansion parameterization}\label{section_local_KL}
Our physical domain $\calD$ can be represented by a finite number, $M$,  of subdomains, i.e., $\ccalD = \cup_{i=1}^M \ccalDi$,
where $\overline{A}$ denotes the closure of the subset $A$. 
We consider the case where the intersection of two subdomains can only be a connected interface with a positive $(N_D-1)$-dimensional measure or an empty set. 
For a subdomain $\calDi$, the set of its boundaries is denoted by $\partial \calDi$, and the set of its neighboring subdomain indices is denoted by $\mathfrak{N}^{(i)}:=\{j|j\in\{1,\ldots,M\}, j\ne i \text{ and } \partial \calDi\cap \partial \calDj\ne \emptyset \}$.
The boundary set $\partial \calDi$ can be split into two parts: external boundaries $\partial \calDi\cap \partial \calD$,  
and interfaces 
$\tauij:=\partial \calDi\cap \partial \calD^{(j)}$ for $j\in \frakN^{(i)}$. 
Grouping all interface indices associated with all subdomains $\{ \calDi \}_{i=1}^M$, we define $\frakN:=\{(i,j)|i\in \{1,2,\ldots,M\} \text{ and }j\in \frakN^{(i)} \}.$  

We introduce decomposed local operators  $\{\calL^{(i)}:= \calL|_{\calDi}\}_{i=1}^M$, $\{\calB^{(i)}:= \calB|_{\calDi}\}_{i=1}^M$ and local functions $\{f^{(i)}:=f|_{\calDi}\}_{i=1}^M$, $\{h^{(i)}:=h|_{\calDi}\}_{i=1}^M$, which are global operators and functions restricted to each subdomain $\calDi$.
The restriction of the field $\kappa(x,\omega)$ to each subdomain is denoted by 
$\kappa^{(i)}(x,\omega):=\kappa(x,\omega)|_{\calDi}$.
For each $i=1,\ldots,M$ and 
$j\in\frakN^{(i)}$, 
$h^{(i,j)}$ denotes an \textit{interface function} defined on the interface $\tauij$, and in this work it is defined as $\hij(x,\omega):=v(x,\omega)|_{\tauij}$, where $v(x,\omega)$ is the solution of the global problem \eqref{eq:original_pde}.
We emphasize that the interface function $h^{(i,j)}$, being the restriction of the global solution on the interface, is dependent on the random input $\omega$.
Each local problem is then defined as: for $i=1,\ldots,M$, find $v^{(i)}(x,\omega): \calDi\times \Omega\to \bbR$ such that
\begin{subequations}
\begin{align}
	&\calL^{(i)} (x,v^{(i)};\kappa^{(i)}(x,\omega)) = f^{(i)}(x)\,,\quad x\in \calDi\,,\label{eq:sub_problem_main} \\
	&\mathcal{B}^{(i)}(x, v^{(i)};\kappa^{(i)}(x,\omega)) = h^{(i)}(x)\,,\quad x\in \partial \calDi\cap \partial \calD\,, \label{eq:sub_problem_boundary}\\
	& \calB^{(i,j)} (x, v^{(i)}; \kappa^{(i)}(x,\omega)) = \hij(x,\omega)\,,\quad x\in \tauij\,,\quad j\in\frakN^{(i)}\label{eq:sub_problem_interface}\,.
\end{align}
\label{eq:sub_problem}
\end{subequations}
Eq.~\eqref{eq:sub_problem_interface} defines the boundary conditions on interfaces and  $\calB^{(i,j)}$ is an appropriate boundary operator posed on the interface $\tauij$. 
With our definition for interface functions, the local problems are consistent with the global problem, i.e.,
\begin{equation*}
	v(x,\omega) = 
	\begin{cases}
	v^{(1)}(x,\omega)\,,\quad x\in \ccalD^{(1)}\,,\\
	\quad \vdots\\
	v^{(M)}(x,\omega)\,,\quad x\in \ccalD^{(M)}\,.
	\end{cases}
\end{equation*}
Cf.\ \cite{liao2015domain,quarteroni1999domain,chen2015local} for detailed discussions for interface functions and boundary conditions for the interfaces. 
         
For each local  random field $\kappa^{(i)}(x,\omega)$ for $i=1,\ldots,M$, its mean function is denoted by $a_0^{(i)}(x)=a_0(x)|_{\calDi}$, where $a_0(x)$ is the mean function of the global field (see \eqref{eq:kle_inf}).
The eigenvalues and the associated orthonormal eigenfunctions of the covariance function posed on each subdomain $\calDi$ are denoted by $\{\lambda_r^{(i)},\psi_r^{(i)}\}_{r=1}^{\infty}$ with $\lambda_1^{(i)}\geq \lambda_2^{(i)}\geq \ldots$, such that
\begin{equation}
\label{eq:local_eigenpair}
	\int_{\calDi} C(x,y) \psi_r^{(i)}(x) \dif x = \lambda_r^{(i)} \psi_r^{(i)}(y)\,,
	\quad x,y\in \calDi\,,
\end{equation}
and $\int_{\calDi}\psi_r^{(i)}(x)\psi_t^{(j)}(x)\dif x=\delta_{rt}$.
The KL expansion of $\kappa^{(i)}(x,\omega)$ can then be written as 
\begin{equation}
\label{eq:local_kl}
\kappa^{(i)}(x,\omega) = a_0^{(i)}(x) + \sum_{r=1}^{\infty}\sqrt{\lambda_r^{(i)}}\psi_r^{(i)}(x)\xi_r^{(i)}(\omega)\,,
\end{equation}
where $\{\xi^{(i)}_r(\omega)\}^{\infty}_{r=1}$ are uncorrelated random variables. 
Each local random field can be approximated by the  truncated KL expansion, 
\begin{equation}
\label{eq:local_kl_trun}
\kappa^{(i)}(x,\omega) \approx a^{(i)}(x,\xi^{(i)}(\omega))=a_0^{(i)}(x) + \sum_{r=1}^{d^{(i)}} \sqrt{\lambda_r^{(i)}}\psi_r^{(i)}(x) \xi_r^{(i)}(\omega)\,,
\quad x\in \calDi\,,
\end{equation}
where $d^{(i)}$ is the number of KL modes retained and 
$\xi^{(i)}(\omega):=[\xi^{(i)}_1(\omega),\ldots, \xi^{(i)}_{d^{(i)}}(\omega)]^T$ whose element is defined as
\begin{equation}
	\label{eq:local_kl_rv}
	\xi^{(i)}_r(\omega):= \frac{1}{\sqrt{\lambda_r^{(i)}}}
	\int_{\calD^{(i)}} (a^{(i)}(x,\xi^{(i)}(\omega)) - a_0^{(i)}(x)) \psi_r^{(i)}\dif x.
\end{equation}
The error of the truncation depends on the amount of total variance captured, $\delta_i:=\sum_{r=1}^{d^{(i)}}\lambda_r^{(i)}/(|\calDi|\sigma^2)$,
and   $d^{(i)}$ needs to be large enough such that $\delta_i$ is larger than  some given threshold $\tolKL$. 

For $i=1,\ldots,M$, the prior distribution of $\xi^{(i)}$ is denoted by $\pi_0(\xi^{(i)})$ with support $I_{\xi^{(i)}}\subset \bbR^{d^{(i)}}$. 
For each $i=1,\ldots,M$ and $j\in \frakN^{(i)}$, $\gij(x,\xi)$ denotes the interface function $\gij(x,\xi):=u(x,\xi)|_{\tauij}$, where $u(x,\xi)$ is the solution of the parameterized global problem \eqref{eq:parameterized_pde}.
The original local problem \eqref{eq:sub_problem} is rewritten as: 
find $u^{(i)}(x,\xi^{(i)})$ mapping $ \calDi\times I_{\xi^{(i)}}$ to $\bbR$ such that
\begin{subequations}
\begin{align}
	&\calL^{(i)} (x,u^{(i)};a^{(i)}(x,\xi^{(i)})) = f^{(i)}(x)\,,\quad x\in \calDi\,,
	\\
	&\calB^{(i)}(x, u^{(i)};a^{(i)}(x,\xi^{(i)})) = h^{(i)}(x)\,,\quad x\in \partial \calDi\cap \partial \calD\,, 
	\\
	& \calB^{(i,j)} (x, u^{(i)}; a^{(i)}(x,\xi^{(i)})) = \gij(x,\xi)\,,\quad x \in \tauij\textrm{ and } j\in\frakN^{(i)}
	\,.
\end{align}
\label{eq:para_subproblem}
\end{subequations}
Defining the observation operator posed on each local subdomain as $\pC^{(i)}:= \pC|_{\calDi}$,
we denote the local forward model as $F^{(i)}(a^{(i)}(x,\xi^{(i)})):= \pC^{(i)}(u^{(i)}(x,\xi^{(i)}))$ and 
the local observation $\bdobs^{(i)}\in \mathbb{R}^{n^{(i)}}$ is defined as $\bdobs^{(i)} = F^{(i)}(a^{(i)}(x,\xi^{(i)})) + \epsilon^{(i)}$, where the local observation noise $\epsilon^{(i)}\sim \mathcal{N}(0,\obsvar^2 \bm{I}_{n^{(i)}})$. 
%As discussed in detail in \cite{chen2015local},  for a given random field with correlation length $L_a$, the decay rate of the eigenvalues (see  \eqref{eq:global_eigenpair_1} or \eqref{eq:local_eigenpair}) depends the relative correlation length, i.e., $L_{a,E} := L_a /L_{E}$, where $L_{E}$ is the diameter of the physical domain $E$. It is shown that long correlation lengths lead to fast decay of eigenvalues, and vice versa. So, when the  correlation length of the global field $\kappa(x,\omega)$ is small, its parameterization over the global domain $\calD$ can be high-dimensional (i.e., $d$ in \eqref{eq:kle_global} is large). To result in a low-dimensional parameterization, we decompose the physical domain into small subdomains, and the relative correlation length then becomes larger. 
Details of our method to efficiently solve the inverse problem
posed on each subdomain are discussed in section \ref{sec:dd_mcmc}, 
where samples of the posterior distribution of each local input $\xi^{(i)}$ are generated.   
The following part of this section is to discuss the procedure of reconstructing the global field
$\kappa(x,\omega)$ with given realizations of  local inputs $\xi^{(i)}$ for $i=1,\ldots, M$.

\subsection{Reconstructed fields} \label{sec:recon}
Letting $(\calDi)^c:=\calD\setminus \calDi$ for $i=1,\ldots,M$, extensions of local mean functions and local eigenfunctions to the global domain $\calD$ are defined as  
\begin{align}
& \tila_0^{(i)}(x): = \begin{cases}
a^{(i)}(x)\,,\quad x\in \calDi\,,\\
0\,,\quad  x\in (\calDi)^c\,,
\end{cases}
\label{eq:stitch_mean}
\\
&\tilpsi_r^{(i)}(x) := 
\begin{cases}
\psi_r^{(i)}(x) \,,\quad x \in \calDi\,,\\
0\,,\quad x\in (\calDi)^c\,.
\end{cases}
\label{eq:stitch_eigenfunction}
\end{align}
For given realizations of the local inputs, two kinds of reconstructed global fields are 
introduced in this section, which are called the stitched field and the assembled  field respectively in the following.

\begin{Def}[The stitched field]
	\label{eq:stitched_field}
When a realization is given for each local input $\xi^{(i)}$ where $i=1,\ldots, M$, 
the stitched field $\brea(x,\xi)$ where $\xi^T:=[(\xi^{(1)})^T,\ldots,(\xi^{(M)})^T]$
and $x\in \calD$,  is defined through directly collecting the corresponding local fields, i.e.,
\[
\brea(x,\xi):= \sum_{i=1}^M \tila^{(i)}(x,\xi^{(i)}),
\]
where $\tila^{(i)}$ is defined as
\begin{align}
\label{eq:local_field}
\tila^{(i)}(x,\xi^{(i)})& :=\tila_0^{(i)}(x) + \sum_{r=1}^{d^{(i)}}\sqrt{\lambda_r^{(i)}}\tilpsi_r^{(i)}(x)\xi_r^{(i)}.
\end{align}
In \eqref{eq:local_field}, $\xi^{(i)}=[\xi^{(i)}_1,\ldots,\xi^{(i)}_{d^{(i)}}]^T$, $\tila_0^{(i)}(x)$ 
and $\tilpsi_r^{(i)}(x)$ are defined in \eqref{eq:stitch_mean}--\eqref{eq:stitch_eigenfunction}, and $\lambda_r^{(i)}$ is defined in \eqref{eq:local_eigenpair}.
\end{Def}

As the stitched field is a direct collection of local fields, it can lead to discontinuities on interfaces, and 
the corresponding  inference results can be misleading. To result in an efficient representation of the global 
field, we define the following assembled field.

\begin{Def}[The assembled field]
\label{def:assembled_field}
When a realization is given for each local input $\xi^{(i)}$ where $i=1,\ldots,M$, the assembled field $\hata(x,\hatxi)$ where $x\in \calD$ is defined as 
\begin{equation}
\label{eq:assembled_kle}
\hata(x,\hatxi) := a_0(x) + \sum_{t=1}^d \sqrt{\lambda_t} \psi_t(x)\hatxi_t\,, 
\end{equation} 
where $\hatxi:=[\hatxi_1,\ldots,\hatxi_d]^T$, the mean function $a_0(x)$ and the eigenpairs $\{\lambda_t, \psi_t(x)\}_{t=1}^{d}$ follow the same settings in \eqref{eq:kle_inf},
and $\{\hatxi_t\}_{t=1}^{d}$ are defined as
\begin{equation}
	\hatxi_t = \frac{1}{\sqrt{\lambda_t}}\sum_{i=1}^M\sum_{r=1}^{d^{(i)}}\sqrt{\lambda_r^{(i)}}\xi_r^{(i)}\int_{\calDi} \tilpsi_r^{(i)}(x)\psi_t(x)\dif x\,,\quad t=1,\ldots,d\,.
	\label{eq:assemble_rv}
\end{equation}
\end{Def}

It can be seen that, the assembled field is represented by the global eigenfunctions, which avoids extra discontinuities on interfaces introduced by the stitched field. In addition, the following theorem shows that the assembled field is the projection of the stitched field
over the space spanned by the global eigenfunctions. 
Therefore, if the truth field is in the space spanned by the global eigenfunctions,
the  assembled field can typically give a better approximation to the truth field than the stitched field.

\begin{thm}
	\label{them:projection}
Let $\mathbb{V}_1$ denote the space spanned by global eigenfunctions $\{\psi_r(x)\}^d_{r=1}$ of \eqref{eq:global_eigenpair_1}.  
The centralized assembled field $\hata(x,\hatxi)-a_0(x)$ is the projection of the centralized stitched field $\brea(x,\xi) - a_0(x)$ over $\mathbb{V}_1$.
\end{thm}

\begin{pf}
For $i=1,\ldots,M$, letting $L_2(\calDi)$ denote the space of square integrable functions over $\calDi$, 
the inner product of any $\mathcal{F}$ and $\mathcal{G}$ belonging to $L_2(\calDi)$ is denoted by  
$\langle \mathcal{F}, \mathcal{G}\rangle_{\calDi}:= \int_{\calDi}\mathcal{F}(x)\mathcal{G}(x)\dif x$.
Denoting $\mathbb{V}_2:={\rm span}\{\tilpsi_r^{(i)}(x), \textrm{ for } i=1,\ldots,M, r=1,\ldots,d^{(i)}\}$,
it can be seen that 
$\hata(x,\hatxi) - a_0(x)\in \mathbb{V}_1$ and $\brea(x,\xi)-a_0(x) \in \mathbb{V}_2$.
From \eqref{eq:assembled_kle}--\eqref{eq:assemble_rv}, the centralized assembled field can be written as 
\begin{align}
    &\hata(x,\hatxi) - a_0(x)  
     =\sum_{t=1}^d \sqrt{\lambda_t}\psi_t(x) \hatxi_t   = \sum_{t=1}^d  \sqrt{\lambda_t}\psi_t(x) \frac{1}{\sqrt{\lambda_t}}\sum_{i=1}^M\sum_{r=1}^{d^{(i)}}\sqrt{\lambda_r^{(i)}}\xi_r^{(i)}(\omega)\int_{\calDi} \tilpsi_r^{(i)}(x)\psi_t(x)\dif x \nonumber \\
    & \quad = \sum_{t=1}^d \sum_{i=1}^M\sum_{r=1}^{d^{(i)}}\sqrt{\lambda_r^{(i)}}\xi_r^{(i)}\langle\tilpsi_r^{(i)}, \psi_t\rangle_{\calDi} \psi_t(x)  =\sum_{t=1}^d \left\langle
    \sum_{i=1}^M [\tila^{(i)}(x,\xi^{(i)}) - \tila_0^{(i)}(x)],\psi_t(x)\right\rangle_{\calD} \psi_t(x) \nonumber \\
    & \quad =\sum_{t=1}^d \left\langle
    \brea(x,\xi) - a_0(x),\psi_t(x)\right\rangle_{\calD} \psi_t(x). 
    \label{eq:projection}
\end{align}
It can be seen from \eqref{eq:projection} that the centralized assembled field is the projection of the stitched field $[\brea(x,\xi) - a_0(x)]$ over $\mathbb{V}_1$.\qed
\end{pf}
For a global  field $a(x,\xi)$, i.e.,  \eqref{eq:kle_global} with a given realization of $\xi$, 
we have that $a(x,\xi) - \hata(x,\hatxi)\in \mathbb{V}_1$. 
For each basis function $\psi_r(x)$ of $\mathbb{V}_1$ (for $r=1,\ldots,d$), Theorem \ref{them:projection} gives that,
\begin{align*}
& \langle \hata(x,\hatxi) - \brea(x,\xi),\psi_r\rangle_{\calD}  = 
\left\langle \sum_{t=1}^d \langle
\brea(x,\xi) - a_0(x),\psi_t\rangle_{\calD}\psi_t + a_0(x) -  \brea(x,\xi)
,\psi_r\right\rangle_{\calD}\\
& = \sum_{t=1}^d \langle
\brea(x,\xi) - a_0(x),\psi_t\rangle_{\calD} \langle \psi_t
,\psi_r\rangle_{\calD}
- 
\langle \brea(x,\xi) - a_0(x)
,\psi_r\rangle_{\calD} =0 \,.
\end{align*}
So, $\langle \hata(x,\hatxi)- \brea(x,\xi), a(x,\xi) - \hata(x,\hatxi)\rangle_{\calD} =0$.
Then, 
\begin{align*}
& \| a(x,\xi) -\brea(x,\xi) \|_{\calD}^2 = 
\| a(x,\xi) - \hata(x,\hatxi) + \hata(x,\hatxi)-  \brea(x,\xi)  \|_{\calD}^2 \\
& \quad = \| a(x,\xi) - \hata(x,\hatxi) \|_{\calD}^2 + 2 \langle \hata(x,\hatxi)- \brea(x,\xi), a(x,\xi) - \hata(x,\hatxi)\rangle_{\calD} + \|\hata(x,\hatxi)-  \brea(x,\xi)  \|_{\calD}^2 \\
&\quad = \| a(x,\xi) - \hata(x,\hatxi) \|_{\calD}^2 + \|\hata(x,\hatxi)-  \brea(x,\xi)  \|_{\calD}^2. 
& 
\end{align*}
Thus, $\|a(x,\xi) -\hata(x,\hatxi) \|_{\calD} \le \| a(x,\xi) -\brea(x,\xi) \|_{\calD}$, which implies that, if 
the given field $a(x,\xi)$ is the truth field of our inverse problem, 
the approximation $\hata(x,\hatxi)$ is typically more accurate than $\brea(x,\xi)$. 

\section{Domain-Decomposed Markov chain Monte Carlo (DD-MCMC)}
\label{sec:dd_mcmc}
%In this section, we demonstrate our Domain-Decomposed Markov chain Monte Carlo (DD-MCMC) algorithm.
%Our goal is to efficiently generate (approximation) samples of  the posterior distribution of the unknown field $\kappa(x,\omega)$ 
%defined in the original global problem \eqref{eq:original_pde} 
Our goal is to efficiently generate samples of  the posterior distribution of the unknown field $a(x,\xi)$ 
 in the  global problem \eqref{eq:parameterized_pde} 
through solving local problems \eqref{eq:para_subproblem}.  
In this section, we first propose a new adaptive Gaussian process (GP) interface model for each 
local problem, and then present our overall DD-MCMC algorithm.

\subsection{Adaptive Guassian process for interface treatments}
\label{sec:gp_fitting}
To include measurement locations, the set consisting of observed data is denoted by 
$S := \{ (x^{s},\dobs^{s})\}_{s=1}^n$, where  $x^{s}$ is the location of the $s$-th sensor and $\dobs^{s}\in \mathbb{R}$ is the observation collected at $x^{s}$. 
The set consisting of all sensor locations is defined by $\bxobs:= \{ x|(x,\dobs)\in S \}$, 
and the observed data are collected as $\bdobs = [\dobs^{1},\ldots,\dobs^{n}]^T$.
For each subdomain $\calD^{(i)}$ (for $i=1,\ldots,M$), 
the set consisting of local observed data is defined as $S^{(i)}:= \{ (x, \dobs)| (x, \dobs)\in S \textrm{ and } x\in \calDi \cup \partial \calDi \}$, of which the size is denoted by $n^{(i)}$,
and $\bdobs^{(i)}\in \bbR^{n^{(i)}}$ collects observed data contained in $S^{(i)}$.

In each local problem \eqref{eq:para_subproblem}, proper interface functions need to be specified.  
%It is known that Guassian process (GP) is widely used to approximate unknown to functions with given data sets. 
%gp的好处variance can be given. we use gp to model interface functions.
%Our goal is to infer the posterior distribution of the unknown field $a(x, \xi(\omega))$ (see \eqref{eq:kle_global}) using local models \eqref{eq:sub1}--\eqref{eq:sub3}, where interface functions need to be specified.
Based on observed data, we build a Gaussian process (GP) model to approximate each interface function, which is a widely used tool to approximate unknown function \cite{rasmussen2006gaussian}.
From Section \ref{section_local_KL}, 
the interface functions to be specified can be considered as the restrictions of the solution of the global problem
associated with the unknown truth sample of the global parameter $\xi$.
That is, for each interface $\tauij$ where $(i,j)\in \frakN$, 
an unknown target function is defined as $\tilgij(x):=\gij(x,\xi)$.
For each target function, its training set is denoted by  $\Lambdaij=\{(x^{s},\dobs^{s})\}_{s=1}^{\nij}$ whose size is denoted by $\nij := |\Lambdaij|$. 
The set consisting of the sensor locations in $\Lambdaij$ is denoted by $\bxobsij:=\{x|(x,\dobs)\in \Lambdaij \}$, 
and  $\bdobs^{(i,j)}\in \mathbb{R}^{\nij}$ collects all observations in $\Lambdaij$. 
Details for constructing the training set and the GP model for the target functions are presented as follows.

A Gaussian process is a collection of random variables, and any finite combinations of these random variables are joint Gaussian distributions.
In our setting, for each of $x$, $\tilgij(x)$ is considered to be a random variable in a Gaussian process.
Each of the prior GP models is denoted by $\tilgij(x) \sim \mathcal{GP}(\mu(x), k(x, x'))$ 
where $\mu(\cdot)$ is the mean function and $k(\cdot, \cdot)$ is the kernel of the Gaussian process.
The Gaussian process is specified by its mean function and kernel function~\cite{rasmussen2006gaussian}.
In this work, we use the Gaussian kernel, 
i.e., $k(x,y) = \sigma_f^2 \exp(-\|x - y\|_2^2/(2l_f^2)$, where the signal variance $\sigma_f$ and the length-scale $l_f$ are both hyper-parameters of the kernel function. Denoting $\gamma = [\sigma_f, l_f]^T$,
for a given training data set $\Lambdaij$, the hyper-parameters can be determined through minimizing the negative log marginal likelihood $\mathcal{M}(\gamma)$:
\begin{align}
\label{eq:gp_hyper}
\mathcal{M}(\gamma) :
&= -\log p(\Lambdaij|\gamma) \\
& = \frac{1}{2} \log \det(K_{\nij}) + \frac{1}{2}(\bdobs^{(i,j)})^T K_{\nij}^{-1} \bdobs^{(i,j)} 
+\frac{\nij}{2}\log(2\pi)\nonumber\,,
\end{align}
where $K_{\nij}$ is the covariance matrix with entries $[K_{\nij}]_{sr} = k(x^{s}, x^{r})$ for $x^{s}, x^{r}\in \bxobsij$ and  
$s,\,r=1,\ldots, \nij$.
Minimizing $\mathcal{M}(\gamma)$ is a non-convex optimization problem and 
we use the Gaussian processes for machine learning toolbox \cite{rasmussen2010gaussian} to solve it. 
% (GPML)

Once the hyper-parameters are determined, 
the conditional predictive distribution for any  $x\in \tau^{(i,j)}$ is given as 
\begin{equation}
	\label{eq:gp_post}
	\tilgij(x)|\Lambdaij, \gamma \sim
	\mathcal{N}(\mu_{\nij}(x), \sigma_{\nij}(x)).
\end{equation}
In \eqref{eq:gp_post}, $\mathcal{N}$ is a Gaussian distribution with  mean $\mu_{\nij}(x)$ and  variance  $\sigma_{\nij}(x)$ defined as  
\begin{subequations}
\begin{align}
	& \mu_{\nij}(x) = k_\star^T (K_{\nij} + \sigma_{\text{obs}}^2\bm{I}_{\nij})^{-1}\bdobs^{(i,j)}\,, 
    \label{eq:gp_post_mean}
\\
	& \sigma_{\nij}(x) = k(x, x) - k_\star^T (K_{\nij} + \sigma_{\text{obs}}^2\bm{I}_{\nij})^{-1}k_\star\,,
    \label{eq:gp_post_variance}
\end{align}
\label{eq:gp_post_details}
\end{subequations}
where $k_\star\in \mathbb{R}^{\nij}$
and its entries are defined as $(k_\star)_s = k(x, x^{s})$ for $x^{s}\in \bxobsij$ and $s=1,\ldots,\nij$ (see~\cite{rasmussen2006gaussian}). 

It is clear that the GP interface model \eqref{eq:gp_post} is determined by the data set  $\Lambdaij$.
To result in an effective data set for each interface GP model, an active training method is developed as follows.
First, the set $\Lambdaij$ is initialized with an arbitrary element in $S^{(i)}$ or $S^{(j)}$,
a test set $\Delta^{(i,j)} \subset \tau^{(i,j)}$ is constructed,
and an initial GP model \eqref{eq:gp_post} using $\Lambdaij$ is constructed.
Second, %the GP model \eqref{eq:gp_post} is obtained with the current training data set $\Lambdaij$, 
variances of the current GP model are computed for each test point $x\in \Delta^{(i,j)}$  using \eqref{eq:gp_post_variance}, and the test point with the largest variance is denoted by  
%Third, for each test sample $x^{s}$ in $\Delta^{(i,j)}$, we can compute its variance indicator $\sigma_{\nij}(x^{s})$ with the current GP model and locate the sample in the test data set with maximum posterior variance, i.e., 
	\begin{equation}
	\label{eq:max_var}
	\overline{x} :=\argmax_{x\in\Delta^{(i,j)}}\sigma_{\nij}(x)\,. 
	\end{equation}
Third, letting $\|\cdot\|_2$ denote the standard Euclidean norm, 
the location of the observation which is closest to $\overline{x}$ is identified as
\[
x^{\star}:= \argmin_{x\in \bxobs}\|x - \overline{x}\|_2,
\]
and the data pair $(x^{\star},\dobs^{\star})$ is then selected to augment the training data set $\Lambdaij$.
The second and third steps are repeated until the maximum posterior variance $\sigma^{\max}_{\Delta^{(i,j)}}:=\max_{x\in\Delta^{(i,j)}}\sigma_{\nij}(x)$ is less than a given threshold $\delta_{\text{tol}}$.
This active learning procedure is included in our main algorithm in the next section.

With the interface GP models, the local problem discussed in section \ref{section_local_KL} is reformulated 
as: %When the adaptive process converges, we utilize the posterior mean function as the interface function.
%The altered parameterized forward sub-problem with GP-fitted interface is then defined as: 
find $u_{\text{GP}}^{(i)}(x,\xi^{(i)}): \calDi\times I_{\xi^{(i)}}\to \mathbb{R}$ such that  
\begin{subequations}
\begin{align}
	&\mathcal{L}^{(i)} (x,u_{\text{GP}}^{(i)};a^{(i)}(x,\xi^{(i)})) = f(x)\,,\quad x\in \calDi\,, \\
	&\mathcal{B}^{(i)}(x, u_{\text{GP}}^{(i)};a^{(i)}(x,\xi^{(i)})) = h^{(i)}(x)\,,\quad x\in \partial \calDi\cap \partial \calD\,,\\
	& \mathcal{B}^{(i,j)} (x, u_{\text{GP}}^{(i)}; a^{(i)}(x,\xi^{(i)})) =
	\mu_{\nij}(x)\,,\quad x \in \tau^{(i,j)},
\end{align}
\label{eq:sub_problem_gp_boundary}
\end{subequations}
where $ j\in \frakN^{(i)}$ and $\mu_{\nij}(x)$ is the mean function of GP interface models defined in \eqref{eq:gp_post}--\eqref{eq:gp_post_details}.
%encodes the GP-fitting interface condition.
Up to now, the local forward model based on \eqref{eq:sub_problem_gp_boundary} is denoted 
by $F_{\text{GP}}^{(i)}(a^{(i)}(x,\xi^{(i)}) ):=\pC^{(i)}(u_{\text{GP}}^{(i)}(x,\xi^{(i)}))$ where $\pC^{(i)}$ is the local 
observation operator as discussed in section \ref{section_local_KL}. 

\subsection{DD-MCMC Algorithm}
\label{sec:dd_mcmc_alg}
%So far, given a large-scale inverse problem, we can decompose the problem in the spatial sense and represent the unknown field with local KL random variables.  In this section, we introduce our domain-decomposed MCMC algorithm. The whole procedure has three phases. The first phase sets up the problem, including domain decomposition, global and local KL eigenpairs, and observation data settings. Then, for every interface $\tau^{(i,j)}$ where $(i,j)\in \frakN$, the second phase adaptively constructs its GP-fitting model.  In the last phase, we run inversions over each subdomain using the local parameterized forward model constructed in the second phase and obtain the assembled posterior mean field.
To begin with,  we compute the global KL expansion \eqref{eq:kle_global} of the prior field $\kappa(x,\omega)$ introduced in our original problem \eqref{eq:original_pde}, 
where the corresponding eigenvalues and eigenfunctions $\{ \lambda_r,\psi_r \}_{r=1}^d$ are defined through \eqref{eq:global_eigenpair_1}--\eqref{eq:global_eigenpair_2}, and divide the global domain $\calD$ into $M$ non-overlapping local domains $\{ \calDi \}_{i=1}^M$. 
Then our domain-decomposed Markov chain Monte Carlo (DD-MCMC) approach 
has  the following steps to generate posterior samples for each local problem. 
The first step is to set up local problems for each subdomain $\calDi$.
This includes computing local KL expansions \eqref{eq:local_kl_trun}, where the corresponding eigenvalues and eigenfunctions $\{\lambda_r^{(i)},\psi_r^{(i)}\}_{r=1}^{d^{(i)}}$ are defined through \eqref{eq:local_eigenpair}, and constructing local observation data sets  $S^{(i)}:= \{ (x, \dobs)| (x, \dobs)\in S \textrm{ and } x\in \calDi \cup \partial \calDi \}$, where $S$ is the set of all observed data pairs in $\calD$. 
The second step is to construct interface conditions for local problems using the adaptive Gaussian process model 
developed in section \ref{sec:gp_fitting}. Through this step, the GP models \eqref{eq:gp_post} for inference functions are built with 
essential observation data, and the variance indicator \eqref{eq:max_var} guarantees the accuracy of the interface condition.
In the third step, for each local subdomain $\calDi$, 
%letting $\bdobs^{(i)}$ collects all observations contained in $S^{(i)}$, 
local posterior samples
$\{\xi_r^{(i),s},s=1,\ldots,\NC, r=1,\ldots,d^{(i)}\}$
 are generated using Algorithm \ref{alg:MHMCMC} with local forward models 
$\widetilde{F}_{\text{GP}}^{(i)}$ (see \eqref{eq:sub_problem_gp_boundary}) and  local observational data $\bdobs^{(i)}$.

With samples $\{\xi_r^{(i),s},s=1,\ldots,\NC, r=1,\ldots,d^{(i)}\}$ for each local problem $i=1,\ldots,M$,
posterior samples of the global input field $a(x,\xi)$ (see \eqref{eq:parameterized_pde}) can
be constructed using Definition \ref{def:assembled_field},
and each   posterior sample is given as 
%Once the local posterior samples $\PosSet := \{\xi_r^{(i,s)},s=1,\ldots,\NC, r=1,\ldots,d^{(i)}\}$ are generated, 
%{\color{blue}we can then obtain the posterior assembled field samples as follows.
%Divide the sample set $\PosSet$ as $\PosSet=\cup_{s=1}^{\NC}P^{s}$ where $\PosSets:=\{\xi_r^{(i,s)}, r=1,\ldots,d^{(i)}\}$, 
%then for $s=1,\ldots,\NC$, the posterior assembled field sample $\hata(x,\hatxi^{s})$ can be defined with $\PosSets$ (see Definition \ref{def:assembled_field}),
\begin{equation}
\label{eq:post_asse_samp}
\hata\left(x,\hatxi^{s}\right) = a_0(x) + \sum_{t=1}^d \sqrt{\lambda_t}\psi_t(x)\hatxi_t^{s}\,,
\end{equation}
where each $\hatxi_t^{s}$ is defined through \eqref{eq:assemble_rv}:
\begin{equation}
\label{eq:post_asse_xi}
	\hatxi_t^{s} = \frac{1}{\sqrt{\lambda_t}}\sum_{i=1}^M\sum_{r=1}^{d^{(i)}}\sqrt{\lambda_r^{(i)}}\xi_r^{(i),s}\int_{\calDi} \tilpsi_r^{(i)}(x)\psi_t(x)\dif x. %\quad t=1,\ldots,d\,.
\end{equation}
%}

Details of our DD-MCMC strategy are summarized in Algorithm \ref{alg:parallel_MCMC}.
Here, $\Delta^{(i,j)}\subset \tauij$ is a given test set and $\delta_{\text{tol}}$ is a given threshold for the variance 
of GP interface models. %% conclusion;
As discussed in Section {\color{blue}\ref{section_local_KL}}, the number of KL modes retained depends on the relative correlation length. 
As the  relative correlation length posed on subdomains is clearly larger than that for the global domain, 
the input parameters of local problems ($\xi^{(i)}$ in \eqref{eq:local_kl_trun}) has lower dimensionalities than the original input parameter ($\xi$ in \eqref{eq:kle_global}). So, the local posterior samples can be efficiently generated in DD-MCMC. 
With the local samples, Definition \ref{def:assembled_field} gives the assembled fields, and Theorem \ref{them:projection}
guarantees that each centralized assembled field is the projection of the corresponding direct 
centralized stitched field (see Definition \ref{eq:stitched_field}) over the space spanned by the global eigenfunctions.

\begin{algorithm}[!ht]
    \caption{Domain-Decomposed MCMC (DD-MCMC) method}
    \label{alg:parallel_MCMC}
    \begin{algorithmic}[1]
        \Require{Observed data $S=\{ (x^{s},\dobs^{s})\}_{s=1}^n$, global domain $\calD$, the mean function and the covaraince function $C(x,y)$ of a prior field.}
        \State{Compute the global KL eigenpairs $\{\lambda_r,\psi_r\}_{r=1}^d$  using \eqref{eq:global_eigenpair_1}--\eqref{eq:global_eigenpair_2}.}
        \State{Partition the global domain $\calD$ into $M$ non-overlapping local domains $\{\calDi\}_{i=1}^M$ with interfaces $\tau^{(i,j)}$ for $j\in \mathfrak{N}^{(i)}$ (see the settings in Section \ref{section_local_KL} for details).}
        \State{Compute the local KL eigenpairs $\{\lambda_r^{(i)},\psi_r^{(i)}\}_{r=1}^{d^{(i)}}$ for $i=1,\ldots,M$ using \eqref{eq:local_eigenpair}.} 
        \State{Divide the data set $S$ into $\{S^{(i)}\}_{i=1}^M$ where $S^{(i)}:= \{ (x, \dobs)|x\in \bxobs \text{ and } x\in \calDi\cup \partial \calDi \}$.}
        \For{each interface $\tau^{(i,j)}$ where $(i,j)\in \frakN $}{
        \State{
        Initialize the training set $\Lambdaij$ with an arbitrary data point in $S^{(i)}\cup S^{(j)}$.}
    	\State{Construct a finite test set $\Delta^{(i,j)} \subset \tau^{(i,j)}$.}
        %\STATE{Train the GP model $\tilgij(x^{\star})|\Lambdaij,(x^{\star}) \sim\mathcal{N}(\mu_{\nij},\sigma_{\nij})$ for arbitrary $x^{\star}\in \tau^{(i,j)}$ using \eqref{eq:gp_post} {\color{blue} with hyper-parameters trained by \eqref{eq:gp_hyper}}.}
        \State{Build a GP interface model 	$\tilgij(x)|\Lambdaij, \gamma \sim
	\mathcal{N}(\mu_{\nij}(x), \sigma_{\nij}(x))$ (see \eqref{eq:gp_hyper}--\eqref{eq:gp_post}).}
        \State{Obtain the maximum posterior variance $\sigma_{\Delta^{(i,j)}}^{\max}:=\max_{x\in\Delta^{(i,j)}}\sigma_{\nij}(x).$}
 %       \WHILE{
 %       $\sigma_{\Delta^{(i,j)}}^{\max}\ge \delta_{\text{tol}}$ and $\Lambdaij\subset S$}
         \While{$\sigma_{\Delta^{(i,j)}}^{\max}\ge \delta_{\text{tol}}$}
         \State{Find $\overline{x}:= \argmax_{x\in \Delta^{(i,j)}}\sigma_{\nij}(x)$ using \eqref{eq:max_var}.}
         \State{Find
        	$x^{\star}:= \argmin_{x\in \bxobs}\|x - \overline{x}\|.$}
            \State{Update the training set: $\Lambdaij=\Lambdaij\cup \{ (x^\star, \dobs^{\star})\}$, where 
       $\dobs^{\star}$ is the observation collected at $x^{\star}$.}
       \State{Go back to line 8.}
        \EndWhile
        }
        \EndFor
%        \STATE{Construct the local forward models $\widetilde{F}_{\text{GP}}^{(i)}$ with $\tau^{(i,j)}$ (see \eqref{eq:sub_problem_gp_boundary}) for $i=1,\ldots,M$. }
        \State{Construct the local forward models $\widetilde{F}_{\text{GP}}^{(i)}$ for $i=1,\ldots,M$ (see \eqref{eq:sub_problem_gp_boundary}). }
    \For{$i=1,\ldots,M$}
    \State{
    	Obtain local posterior samples $\left\{\xi_r^{(i),s}, {\textrm{ for }} s=1,\ldots,\NC, r=1,\ldots,d^{(i)}\right\}$ using Algorithm \ref{alg:MHMCMC} with local model $\widetilde{F}_{\text{GP}}^{(i)}$ and local observation data $\bdobs^{(i)}$.}
    \EndFor
    \State{
  Construct samples of the assembled field $\left\{\hata(x,\hatxi^{s})\right\}_{s=1}^{\NC}$ using \eqref{eq:post_asse_samp}--\eqref{eq:post_asse_xi}. %(see details in Section \ref{sec:dd_mcmc_alg}).
  	%With local posterior mean estimates defined in \eqref{eq:post_mean_estimates},
  	%obtain the assembled posterior mean field $\checkE[\widehat{a}(x,\widehat{\xi})]$ using \eqref{eq:assem_post_mean} (see details in Section \ref{sec:dd_mcmc_alg}).
  }
  \Ensure{Posterior samples $\left\{\hata(x,\hatxi^{s})\right\}_{s=1}^{\NC}$.}
  	%The assembled posterior mean field $\checkE[\widehat{a}(x,\widehat{\xi})]$.
    \end{algorithmic}
\end{algorithm}
 
\section{Numerical study} 
\label{sec:numerical}
In this section, numerical experiments are conducted to illustrate the effectiveness of our  domain-decomposed Markov chain Monte Carlo (DD-MCMC) approach. 
%In this section, we present and discuss the proposed domain-decomposed Markov chain Monte Carlo method on a test problem involving an elliptic PDE with random coefficients.
Setups for our three test problems are described in Section \ref{sec:num_setup}.
Effects of the Gaussian process (GP) interface treatments are discussed  in Section \ref{sec:gp_fitting_num},
and the overall inference results of DD-MCMC are discussed in Section \ref{sec:3_sub}. 

% Detailed Problem setup
\subsection{Setup for test problems}
\label{sec:num_setup}

The numerical examples considered are steady flows in porous media. 
Letting $a(x,\xi)$ denote an unknown permeability field and $u(x,\xi)$ the pressure head, we consider the following diffusion equation, 
\begin{align}
-\nabla\cdot (a(\bx,\bxi)\nabla u(\bx,\xi)) &= f(\bx)\,,\quad \bx\in \calD\,. \label{eq_diff}
\end{align}
Here, the physical domain  considered is $\calD=(0,3)\times (0,1)\subset \bbR^2$,
and the homogeneous Dirichlet boundary condition is specified on the left and right boundaries and the homogeneous Neumann boundary condition is specified on the top and bottom boundaries, i.e., 
% boundary conditions 
\begin{align*}
& u(\bx,\bxi) = 0\,,\quad \bx\in \{0\}\times[0,1]\\
& u(\bx,\bxi) = 0\,,\quad \bx\in \{3\}\times[0,1]\,,\\
& a(\bx,\bxi)\nabla u(\bx,\bxi) \cdot \bm{n}(\bx) =0 \,,
\quad \bx\in \{(0,3)\times \{0\}\}\cup \{(0,3)\times \{1\} \} \,,
\end{align*}
where  $\bm{n}(\bx)$ is the unit normal vector to the %Neumann 
boundary. 
The source term is specified as   
\[
f(\bx) = 3\exp\left(-\|\bx^{sr} - \bx\|_2^2\right)\,,
\]
where $\bx^{sr} = [\bx^{sr}_1,\bx_2^{sr}]^T$ denotes the center of contaminant and it is set to $\bx^{sr} = [1.5,0.5]^T$. 
%\textcolor{red}{In this work, the center of contaminant is set to the center of the considered domain, i.e., $\bx^{sr} = (1.5,0.5)$.
Given a realization of $\bxi$, the bilinear finite element method \cite{elman14finite, ifiss} is applied to solve this diffusion equation, where 
a uniform $97\times 33$ grid (the number of the degrees of freedom is 3201) is used.
Our deterministic global forward model $F(\xi)$ is defined to be a set collecting solution values corresponding to measurement sensors, which are uniformly located in the tensor product $\{x_1^{i}\}\otimes \{x_2^{j}\}$ of the one-dimensional grids: 
$x_1^{i}=0.125i, i=1,\ldots,23$, $x_2^{j}=0.125j, j=1,\ldots,7$, where 161 sensors in total are included. The measurement noises are set to independent and identically distributed Gaussian distributions with mean zero, and  the standard deviation is set to $1\%$ of the mean observed value.

\begin{figure}[!htp]
    \centerline{
    \begin{tabular}{c}
    \includegraphics[width=0.80\textwidth]{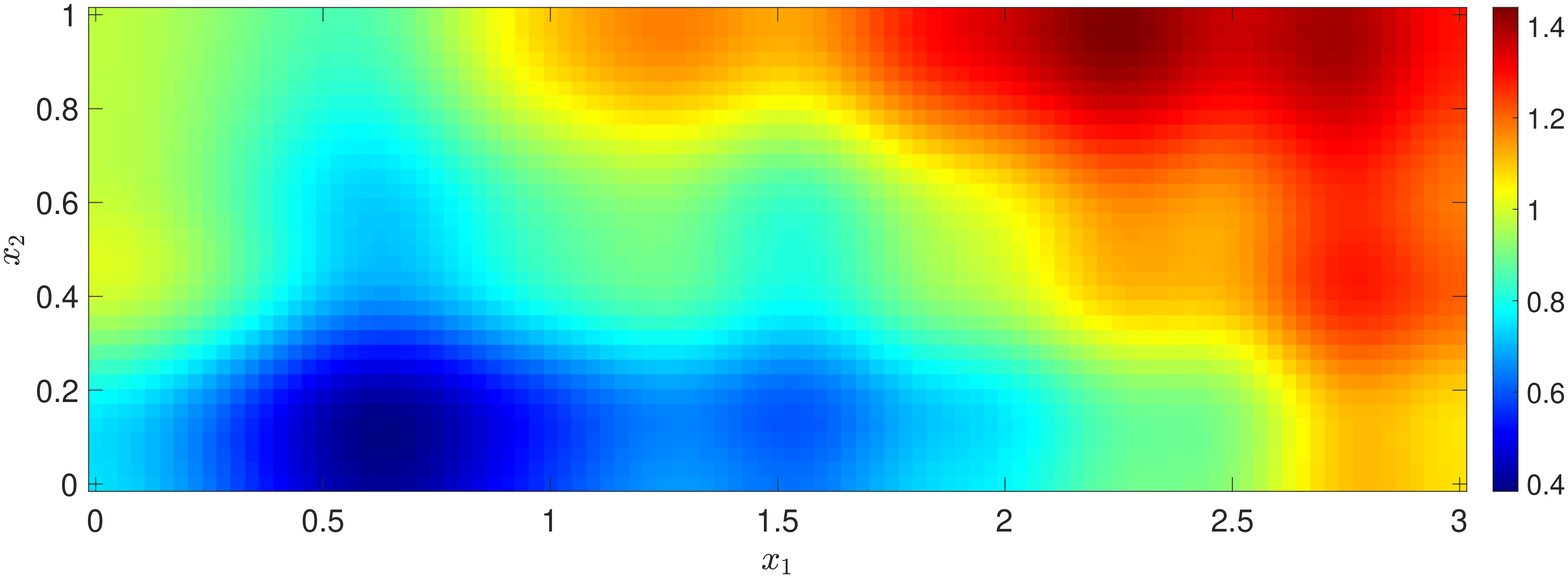} \\
    (a) The truth permeability \\
    \includegraphics[width=0.80\textwidth]{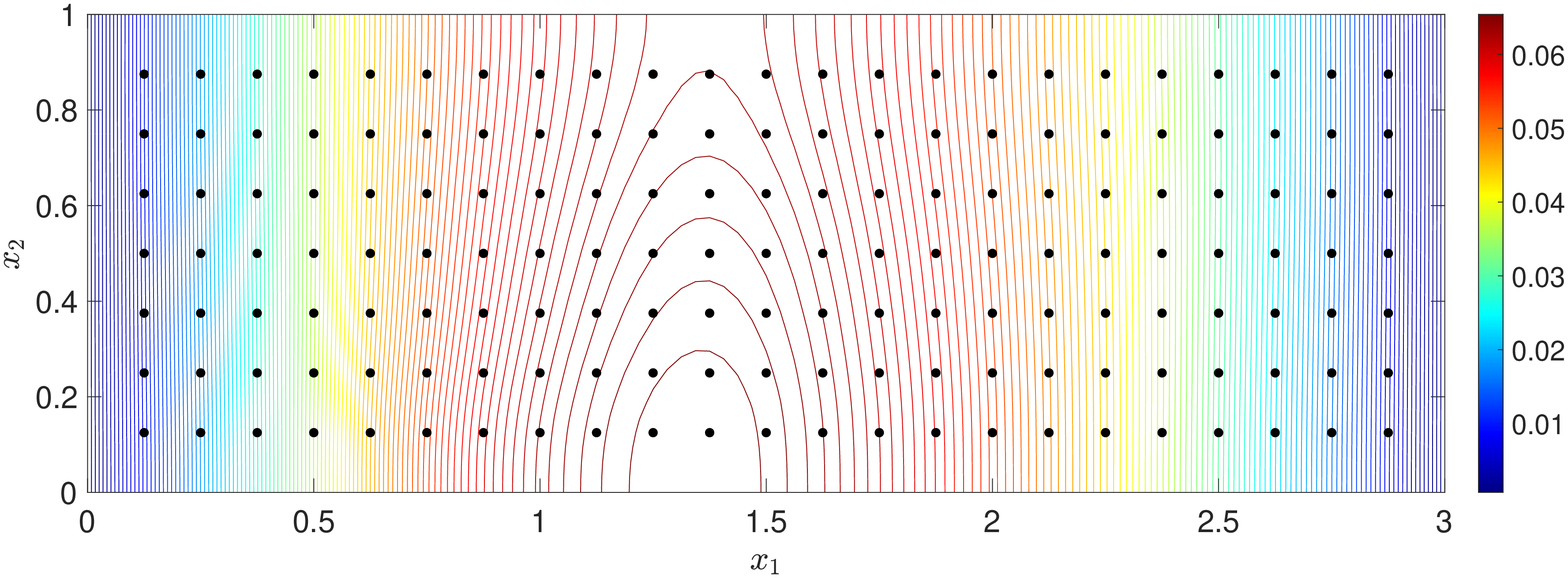}\\
     (b) The pressure field and sensors
    \end{tabular}}
    \caption{Test problem one setup ($\covl=2$).}
    \label{fig:setup_2}
\end{figure}

\begin{figure}[!htp]
    \centerline{
    \begin{tabular}{c}
    \includegraphics[width=0.80\textwidth]{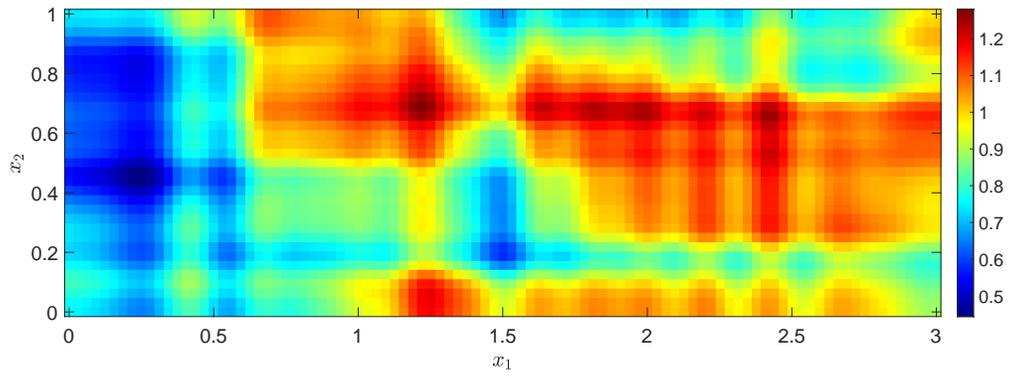} \\
    (a) The truth permeability \\
    \includegraphics[width=0.80\textwidth]{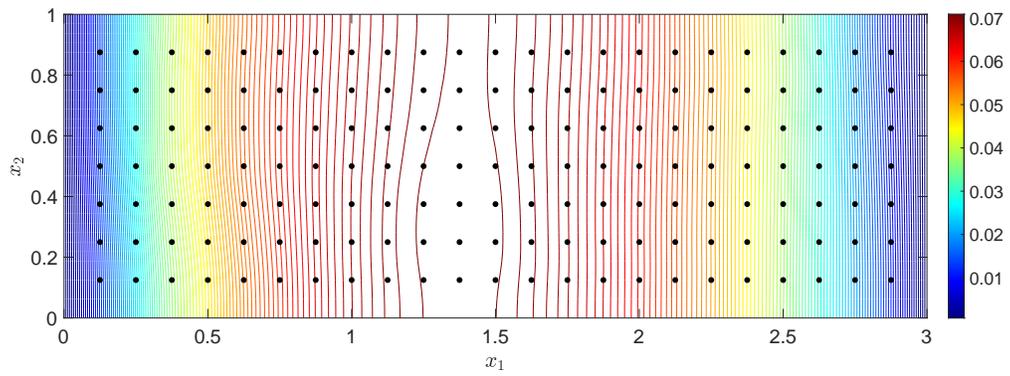}\\
     (b) The pressure field and sensors
    \end{tabular}}
    \caption{Test problem two setup ($\covl=1$).}
    \label{fig:setup_1}
\end{figure}

\begin{figure}[!htp]
    \centerline{
    \begin{tabular}{c}
    \includegraphics[width=0.80\textwidth]{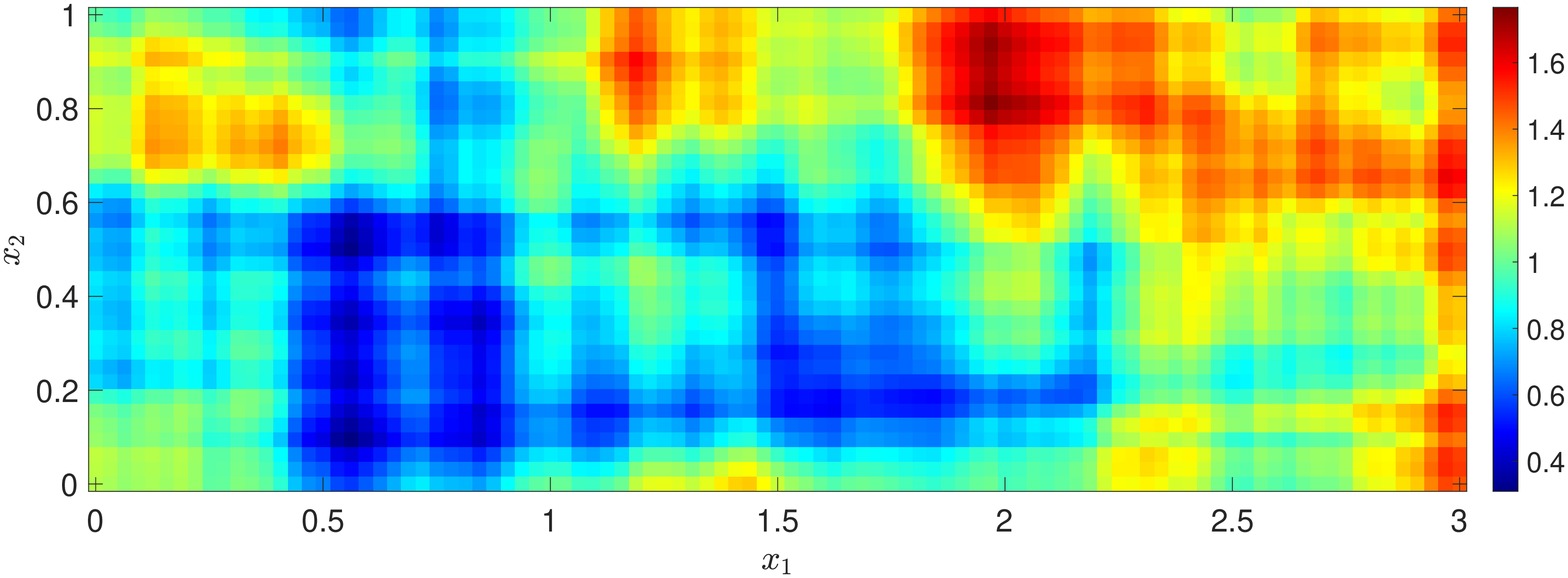} \\
    (a) The truth permeability \\
    \includegraphics[width=0.80\textwidth]{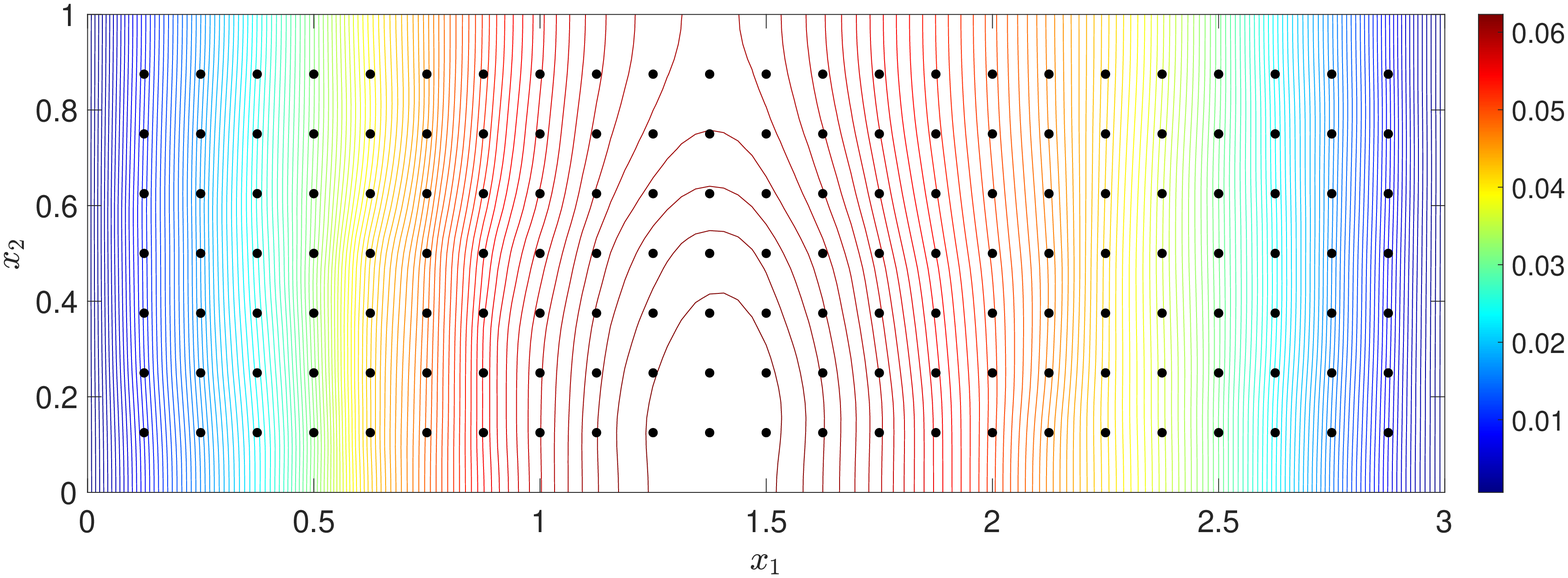}\\
     (b) The pressure field and sensors
    \end{tabular}}
    \caption{Test problem three setup ($\covl=0.5$).}
    \label{fig:setup_05}
\end{figure}

In \eqref{eq_diff}, we set the permeability field $a(x,\xi)$ to a truncated KL expansion of a random field with mean function $a_0(x)$, standard deviation $\sigma$ and covariance function
\begin{eqnarray*}
 C(x,y)=\sigma^2 \exp\left(-\frac{|x_1-y_1|}{L}-\frac{|x_2-y_2|}{L}\right), \label{eq:covariance}
\end{eqnarray*}
where $L$ is the correlation length,
and set $a_0(x)=1$ and $\sigma=0.25$ in the following numerical studies. 
The priori distributions of $\{\xi_r\}^d_{r=1}$ (see \eqref{eq:kle_global}) are set to be independent uniform distributions with range $I=[-1,1]$, and the support of $\xi$ is then $I_{\xi} = I^d$. 
As usual, we set $d$ large enough, such that $\tolKL=95\%$ (see \eqref{eq_tolKL}) of the total variance of the covariance function are captured.  Three test problems are considered in this section, which are associated with  
three different values of the correlation length $\covl$: $2$, $1$ and $0.5$,
and the number of global KL terms retrained are $d=27$, $d=87$ and $d=307$ respectively. 
Figure \ref{fig:setup_2}, Figure \ref{fig:setup_1} and Figure \ref{fig:setup_05} show the truth permeability fields and sensor locations with the corresponding pressure fields for the three test problems respectively. 

For each subdomain $\calD^{(i)}$ ($i=1,2,3$), the local KL expansion  is computed (see \eqref{eq:local_kl_trun}).
To capture $95\%$ of the total variance for each subdomain, the number of local KL modes retained is $11$ for $L=2$, that is 33 for $L=1$ and that is $109$ for $L=0.5$ (as the subdomains have the same dimensionality, the numbers of their corresponding KL modes retained are the same for a given correlation length).  
The priori distributions of the local parameters $\{\xi_r^{(i)}\}_{r=1}^{d^{(i)}}$ (see \eqref{eq:local_kl}) for $i=1,2,3$ are set to be independent uniform distributions with range $I=[-1,1]$, and the support of $\xi^{(i)}$ is then $I_{\xi^{(i)}} = I^{d^{(i)}}$. For each local subdomain, the local problem \eqref{eq:sub_problem_gp_boundary} is discretized with the bilinear finite element method with a uniform  $33\times 33$ grid (the number of the spatial degrees of freedom is 1089).  All results of this paper are obtained in MATLAB on a workstation with 2.20 GHz Intel(R) Xeon(R) E5-2630 CPU.
As solving the linear system associated with the global model \eqref{eq:parameterized_pde} takes around $6.256\times 10^{-2}$ %0.06256 
seconds and solving that associated with the local model \eqref{eq:sub_problem_gp_boundary} is around  $3.845\times 10^{-3}$ %0.003845 
seconds, we define the computational cost to conduct a local forward model evaluation as one cost unit, and consider the cost for each global model evaluation to be 16.25 cost units.

For comparison, the standard MCMC method (Algorithm \ref{alg:MHMCMC}) is applied with the global forward model \eqref{eq:parameterized_pde}, which is referred to as the global MCMC (G-MCMC) in the following. For both G-MCMC and our DD-MCMC (Algorithm \ref{alg:parallel_MCMC}), the proposal distribution (line 3 of Algorithm \ref{alg:MHMCMC}) is set to the symmetric Gaussian distribution, i.e., $Q(\xi^{\star}|\xi^{s}) = \mathcal{N}(\xi^{s}, \beta \bm{I})$ where $\bm{I}$ is an identity matrix and $\beta$ is the stepsize, which is specified  for each test problem in Section \ref{sec:3_sub}. For test problem one ($L=2$), the number of posterior samples $\NC$ generated by DD-MCMC is set to $1\times 10^4$ (see line 19 of Algorithm \ref{alg:parallel_MCMC}); for   test problem two ($L=1$), that is set to $2\times 10^4$;  for test problem three ($L=0.5$), that is set to $4\times 10^4$. For a fair comparison,  the number of posterior samples generated by G-MCMC is set to $1\times 10^3$, $2\times 10^3$ and $4\times 10^3$ for test problem one, test problem two and test problem three respectively, such that the costs of DD-MCMC and G-MCMC are approximately equal.

\subsection{Results for the interface treatment}
\label{sec:gp_fitting_num}
To construct GP models for the interface functions (discussed in Section \ref{sec:gp_fitting}), the test sets $\Delta^{(1,2)}$ and $\Delta^{(2,3)}$ (see line 7 of Algorithm \ref{alg:parallel_MCMC}) are set to the grid points on the interfaces, where each interface has $33$ grid points,
and the threshold $\delta_{\text{tol}}$ for the maximum posterior variance (see line 10 of Algorithm \ref{alg:parallel_MCMC}) is set to $10^{-7}$. 
The maximum numbers of training data points required for the active learning procedure (i.e., $|\Lambdaij|$ for the last iteration step in line 13 of Algorithm \ref{alg:parallel_MCMC}) and the corresponding maximum posterior variance (i.e., $\sigma_{\Delta^{(i,j)}}^{\max}$ in line 9 of Algorithm \ref{alg:parallel_MCMC} after the last iteration step)
are shown in Table \ref{tab:gp_fitting_N_max}, where it can be seen that these numbers are small---at most five training data points are required to reach the desired threshold for the three test problems.  

To assess accuracy of the interface treatment, we compute the difference between the GP interface models and the exact interface functions associated with the truth permeability fields of the three test problems. For each interface $\tau^{(i,j)}$,
the relative interface error is computed through 
\begin{eqnarray*}
\reint^{(i,j)} := \| g^{(i,j)}(x,\xi) - \widetilde{g}^{(i,j)}(x) \|_2/\|g^{(i,j)}(x,\xi)\|_2, 
\end{eqnarray*}
where $g^{(i,j)}(x,\xi)$ is the exact interface function defined as $g^{(i,j)}(x,\xi):=u
(x,\xi)|_{\tau^{(i,j)}}$ and $\widetilde{g}^{(i,j)}(x)$ is the trained GP interface model (see line 8 of Algorithm \ref{alg:parallel_MCMC}), and the parameter value $\xi$ is associated with each of the truth field. 
Moreover, for each local  subdomain $\calD^{(i)}$ ($i=1,2,3$), 
the relative state errors of local solutions obtained with the GP interface models are also assessed, which are computed through 
\[
	\restate^{(i)} := \| u_{\text{GP}}^{(i)}(x,\xi^{(i)}) - u^{(i)}(x,\xi^{(i)}) \|_2/\|u^{(i)}(x,\xi^{(i)})\|_2\,,
\]
where $u_{\text{GP}}^{(i)}(x,\xi^{(i)})$ is the local solution defined in \eqref{eq:sub_problem_gp_boundary}, $u^{(i)}(x,\xi^{(i)})$ is the exact local solution which is defined in \eqref{eq:para_subproblem}, and $\xi^{(i)}$ is the local KL random variable (see \eqref{eq:local_kl_rv}) associated with $\xi$.
Table \ref{tab:gp_error} shows the relative interface errors and state errors, where it can be seen that these errors are all very small. 

%Interface relative errors of different interfaces and the state relative errors over different subdomains for different $L$ are displayed in Table \ref{tab:gp_error}.
%With a larger correlation length, the fitting performance slightly increases due to the decrease in model complexity.  Overall, the model error introduced by the GP-fitting interface is beyond the noise level of the measurement noise.

\begin{table}[!htp]
	\caption{Maximum number of training data points and the correspoding maximum posterior variances for the three test problems. }
	\centering
	\begin{tabular}{c|c c c}
		\hline
		 $L$& $2$  &$1$  & $0.5$ \\
		\hline
         $|\Lambda^{(1,2)}|$ & 4 & 4& 4 \\  
        $|\Lambda^{(2,3)}|$ & 4 & 4 & 5\\
		\hline
		$\sigma_{\Delta^{(1,2)}}^{\max}$ & $ 3.979\times 10^{-13}$  & $ 2.700\times 10^{-13}$& $ 4.105\times 10^{-9}$\\  
        $\sigma_{\Delta^{(2,3)}}^{\max}$ & $ 1.856\times 10^{-10}$ & $ 1.160\times 10^{-10}$ & $ 7.555\times 10^{-9}$\\
		\hline
	\end{tabular}
	\label{tab:gp_fitting_N_max}
\end{table}

\begin{table}[!htp]
	\caption{Relative errors for different interfaces}
	\centering
	\begin{tabular}{c|c|ccc}
		\hline
		 &$L$& $2$  &$1$  & $0.5$ \\
		\hline
		\multirow{2}*{$\reint$} &
        $\tau^{(1,2)}$ & $3.112 \times 10^{-3}$ & $3.825\times 10^{-3}$& $4.628 \times 10^{-3}$ \\ 
		~ &  $\tau^{(2,3)}$ & $1.580 \times 10^{-3}$ & $1.925\times 10^{-3}$& $2.345 \times 10^{-3}$ \\  
\hline
\multirow{2}*{$\restate$} & $\calD^{(1)}$ & $ 1.631\times 10^{-3}$ & $ 1.728 \times 10^{-3}$ & $ 2.385 \times 10^{-3}$\\
~ & $\calD^{(2)}$ & $2.743 \times 10^{-4}$ & $3.056 \times 10^{-4}$ & $4.455\times 10^{-4}$\\
 ~ & $\calD^{(3)}$ & $8.375 \times 10^{-5}$ & $1.083 \times 10^{-4}$ & $1.344 \times 10^{-4}$\\
\hline
	\end{tabular}
	\label{tab:gp_error}
\end{table}
\subsection{Performance of DD-MCMC}
\label{sec:3_sub}
%In this section, we present and discuss the performance of the proposed DD-MCMC algorithm.
For the three test problems, values of the stepsize ($\beta$ is introduced in Section \ref{sec:num_setup}) are set as follows, such that the acceptance rate of DD-MCMC and G-MCMC is appropriate,  which  is defined by the number of accepted samples (line 7 of Algorithm \ref{alg:MHMCMC}) divided by the total sample size. 
For test problem one ($L=2$), the stepsize for DD-MCMC is set to $\beta =0.05$, and that for G-MCMC is set to $\beta = 0.07$. For test problem two ($L=1$), that is set to $\beta = 0.05$ for both DD-MCMC and G-MCMC. For test problem three ($L=0.5$), that is set to $\beta = 0.05$ for DD-MCMC and $\beta = 0.04$ for G-MCMC. Acceptance rates for the three test problems are shown in Table \ref{tab:ar}, which are consistent with the settings discussed in \cite{roberts2001optimal}. 

\begin{table}[!htp]
	\caption{Acceptance rates for the three test problems. }
	\centering
	\begin{tabular}{c|ccc|c}
		\hline
		\multirow{2}{*}{$L$} & \multicolumn{3}{c|}{DD-MCMC} & \multirow{2}{*}{G-MCMC} \\ %\cline{2-4}
		& $\calD^{(1)}$& $\calD^{(2)}$ & $\calD^{(3)}$ & \\                      
		\hline
		2& 48.73\% & 16.63\% & 15.97\%  & 17.80\% \\
		1& 42.10\% & 24.68\% & 12.51\%  & 15.00\% \\
		0.5& 47.73\% & 23.79\% & 19.37\%  & 12.60\% \\
		\hline
	\end{tabular}
	\label{tab:ar}
\end{table}

Once samples of the posterior distributions are obtained, we compute the mean and variance estimates of the permeability fields as follows. For our DD-MCMC (Algorithm \ref{alg:parallel_MCMC}), the outputs are denoted by $\{\hata(x,\hatxi^{s})\}_{s=1}^{\NC}$ which are the posterior assembled fields (see \eqref{eq:post_asse_samp}--\eqref{eq:post_asse_xi}) to approximate the unknown permeability field. The corresponding mean and variance estimates are computed through 
\begin{align}
&  
\checkE\left[\hata(x,\hatxi)\right]
:=\frac{1}{\NC}\sum_{s=1}^{\NC}\left[\hata(x,\hatxi^{s})\right], \label{eq:assem_post_mean}\\
&  \checkV\left[\hata(x,\hatxi)\right]:= \frac{1}{N}\sum_{s=1}^N \left[ \hata(x,\hatxi^{s}) - \checkE[\hata(x,\hatxi)] \right]^2.
\label{eq:assemble_var}
\end{align}
For the local posterior samples $\{\bxi^{(i),s},s=1,\ldots,\NC\}$ for $i=1,\ldots, M$ (corresponding to line 19 of Algorithm \ref{alg:parallel_MCMC}), the mean estimate of the local permeability field is obtained by putting samples of each local field (i.e.\ $\{a^{(i)}(x,\xi^{(i),s}), s=1,\ldots,\NC\}$ defined in \eqref{eq:local_kl_trun}) into \eqref{eq:assem_post_mean}. 
%the local posterior mean fields are approximated as 
%\begin{align}
%\checkE[a^{(i)}(x,\xi^{(i)})]  := 
%\frac{1}{N}\sum_{s=1}^N a^{(i)}(x,\xi^{(i),s})\,,\quad i=1,\ldots, M\,.\label{eq:local_post_mean} 
%\end{align}
Denoting the posterior samples of the stitched field as $\brea(x,\xi^{s}):=\sum_{i=1}^M \tila^{(i)}(x,\xi^{(i),s})$, 
where $\tila^{(i)}(x,\xi^{(i),s})$ is the extension of the local field $a^{(i)}(x,\xi^{(i),s})$ following \eqref{eq:stitch_mean}--\eqref{eq:stitch_eigenfunction}. The global mean and variance estimates based on the stitched field, denoted by $\checkE\left[\breve{a}(x,\xi)\right]$ and $\checkV\left[\breve{a}(x,\xi)\right]$ respectively,  are obtained by putting the samples $\{\tila^{(i)}(x,\xi^{(i),s}),s=1,\ldots,\NC\}$ into \eqref{eq:assem_post_mean}--\eqref{eq:assemble_var}. Moreover, 
%the stitched posterior mean field can then be estimated via 
%\begin{align}
%\checkE[\brea(x,\xi)]:=
%\frac{1}{N}\sum_{s=1}^N \brea(x,\xi^{s})\,.\label{eq:stitch_post_mean}
%\end{align}
%Moreover, the corresponding assembled posterior samples $\{\hatxi^s\}_{s=1}^N$ can be computed via \eqref{eq:post_asse_samp}--\eqref{eq:post_asse_xi}.
%Similarly, the assembled posterior mean field and the assembled posterior variance field are approximated with
%\begin{align}
%&  
%\checkE[\hata(x,\hatxi)]
%:=\frac{1}{\NC}\sum_{s=1}^{\NC}[\hata(x,\hatxi^{s})]\,, %\label{eq:assem_post_mean}\\
%&  \checkV[\hata(x,\hatxi)]:= \frac{1}{N}\sum_{s=1}^N \left[ \hata(x,\hatxi^{s}) - \checkE[\hata(x,\hatxi)] \right]^2\,,
%%\label{eq:assemble_var}
%\end{align}
%and the assembled posterior standard deviation field can be defined via $\check{\sigma}[\hata(x,\hatxi)]:=\sqrt{\checkV[\hata(x,\hatxi)]}$.
for samples $\{\bxi^{s},s=1,\ldots,\NC\}$ generated by G-MCMC, the corresponding samples of the global permeability field
are denoted by $\{a(x,\xi^{s}), s=1,\ldots,\NC\}$ (see \eqref{eq:kle_global}). The global mean and variance estimates, denoted by $\checkE\left[a(x,\xi)\right]$ and $\checkV\left[a(x,\xi)\right]$ respectively, are  assessed through putting $\{a(x,\xi^{s}), s=1,\ldots,\NC\}$ into \eqref{eq:assem_post_mean}--\eqref{eq:assemble_var}.
%We define the following posterior quantities to access the inversion results. In the global sense, for a set of posterior samples  $\{\bxi^{s},s=1,\ldots,\NC\}$ generated by the G-MCMC algorithm, the global posterior mean field and the global posterior variance field are estimated by 
%\begin{align}
%& \checkE[a(x,\xi)]:= 
%\frac{1}{N}\sum_{s=1}^N a(x,\xi^{s}) \,,\label{eq:global_post_mean}\\
%& \checkV[a(x,\xi)]:= \frac{1}{N}\sum_{s=1}^N \left[a(x,\xi^s) - \checkE[a(x,\xi)]\right]^2\,, 
%\label{eq:global_post_var}
%\end{align}

\begin{figure}[!htp]
    \centerline{
    \begin{tabular}{c}
    \includegraphics[width=0.80\textwidth]{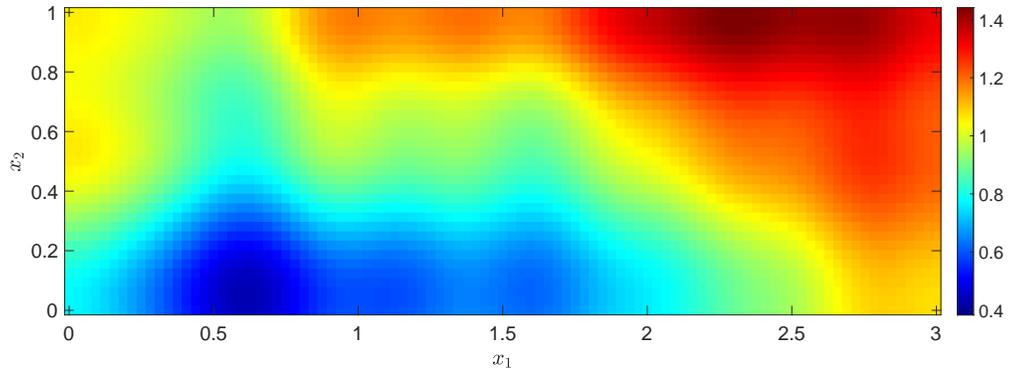} \\
    (a) Mean $\checkE\left[\hata(x,\hatxi)\right]$, DD-MCMC (assembled).\\
	\vspace{0.05mm} \\
    \includegraphics[width=0.80\textwidth]{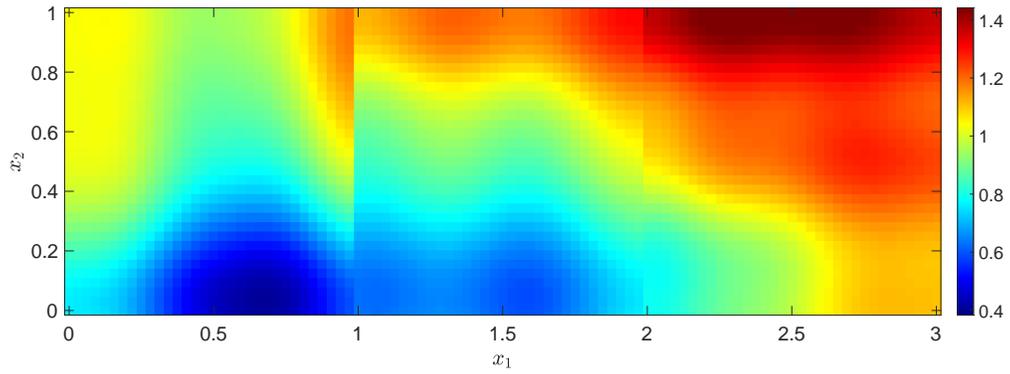}\\
     (b) Mean $\checkE\left[\breve{a}(x,\xi)\right]$, DD-MCMC (stitched).\\
	 \vspace{0.05mm} \\
	 \includegraphics[width=0.80\textwidth]{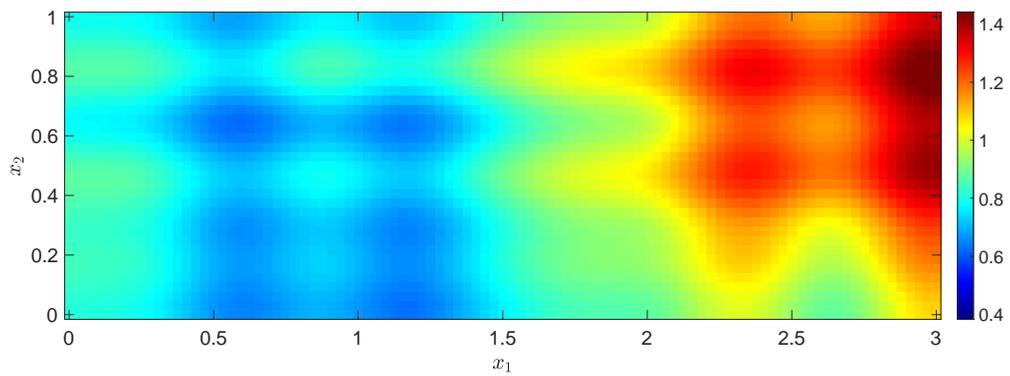}\\
     (c) Mean $\checkE\left[ a(x,\xi) \right]$, G-MCMC.
    \end{tabular}}
    \caption{Estimated mean fields for test problem one ($L=2$).}
    \label{fig:mean_compare_L_2}
\end{figure}

\begin{figure}[!htp]
    \centerline{
    \begin{tabular}{c}
    \includegraphics[width=0.80\textwidth]{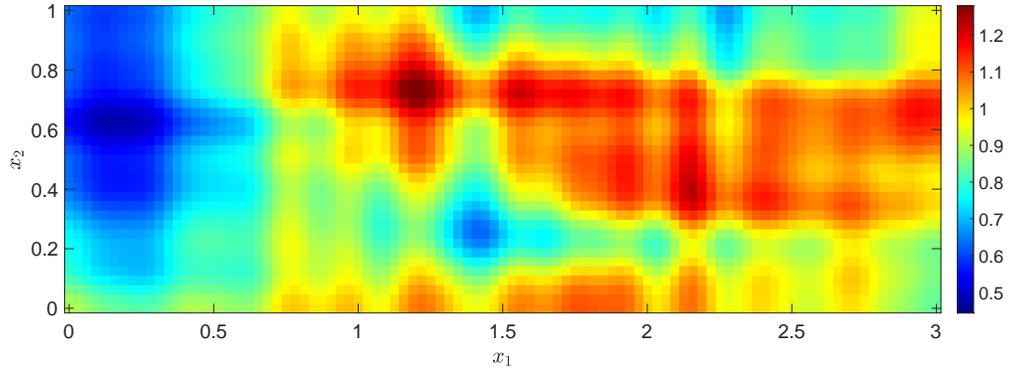} \\
    (a) Mean $\checkE\left[\hata(x,\hatxi)\right]$, DD-MCMC (assembled).\\
	\vspace{0.05mm} \\
    \includegraphics[width=0.80\textwidth]{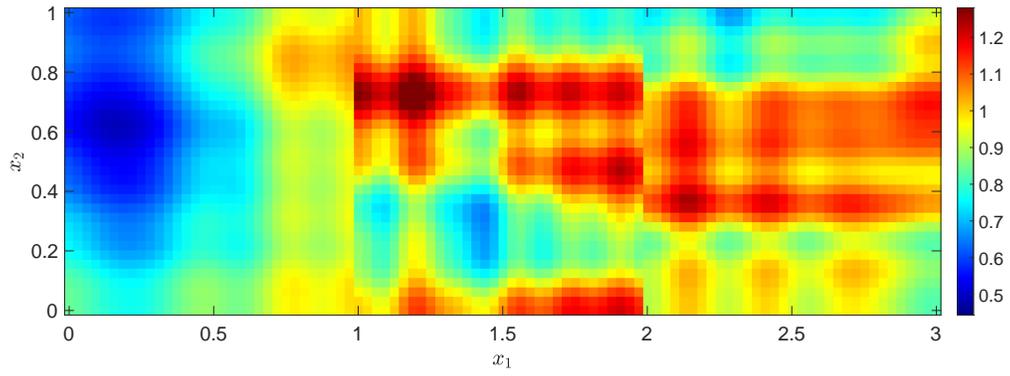}\\
     (b) Mean $\checkE\left[\breve{a}(x,\xi)\right]$, DD-MCMC (stitched).\\
	 \vspace{0.05mm} \\
	 \includegraphics[width=0.80\textwidth]{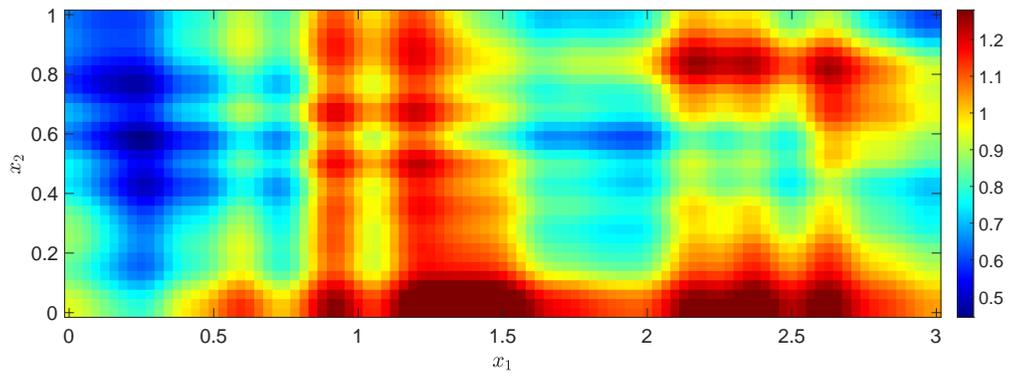}\\
     (c) Mean $\checkE\left[ a(x,\xi) \right]$, G-MCMC .
    \end{tabular}}
    \caption{Estimated mean fields for test problem two ($L=1$).}
    \label{fig:mean_compare_L_1}
\end{figure}

\begin{figure}[!htp]
    \centerline{
    \begin{tabular}{c}
    \includegraphics[width=0.80\textwidth]{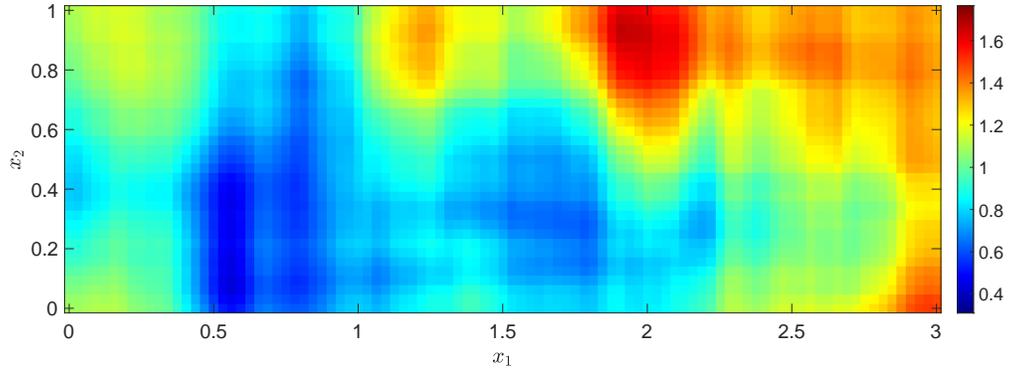} \\
    (a) Mean $\checkE\left[\hata(x,\hatxi)\right]$, DD-MCMC (assembled).\\
	\vspace{0.05mm} \\
    \includegraphics[width=0.80\textwidth]{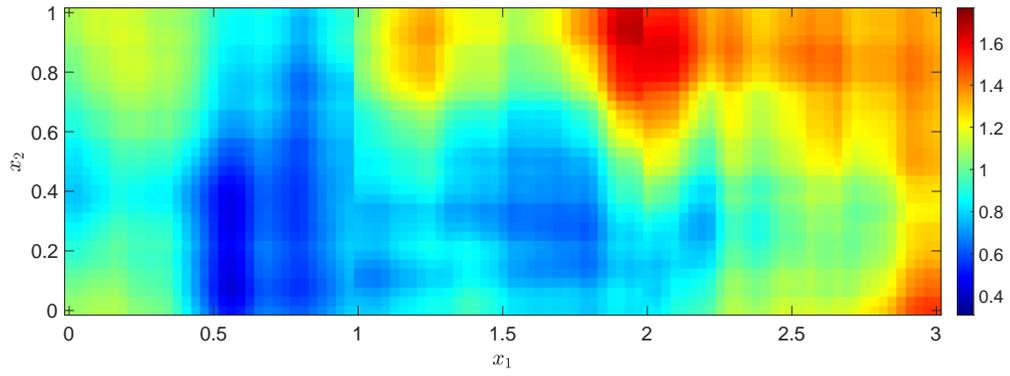}\\
     (b) Mean $\checkE\left[\breve{a}(x,\xi)\right]$, DD-MCMC (stitched).\\
	 \vspace{0.05mm} \\
	 \includegraphics[width=0.80\textwidth]{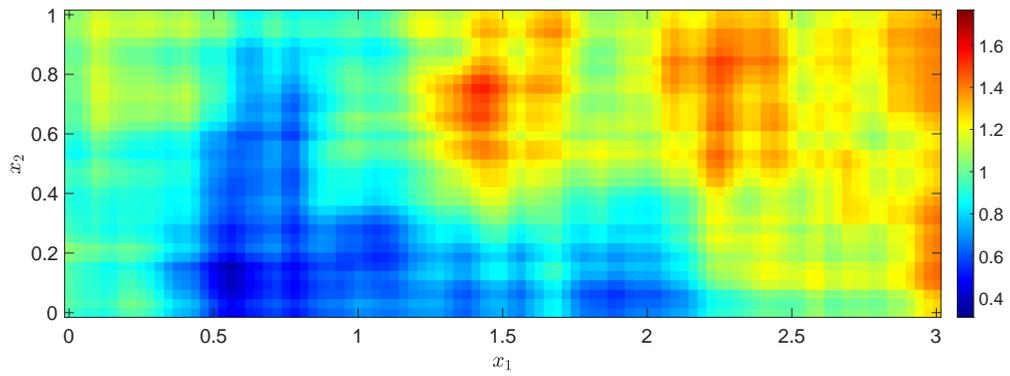}\\
     (c) Mean $\checkE\left[ a(x,\xi) \right]$, G-MCMC.
    \end{tabular}}
    \caption{Estimated mean fields for test problem three ($L=0.5$).}
    \label{fig:mean_compare_L_05}
\end{figure}

\begin{figure}[!htp]
    \centerline{
    \begin{tabular}{c}
    \includegraphics[width=0.80\textwidth]{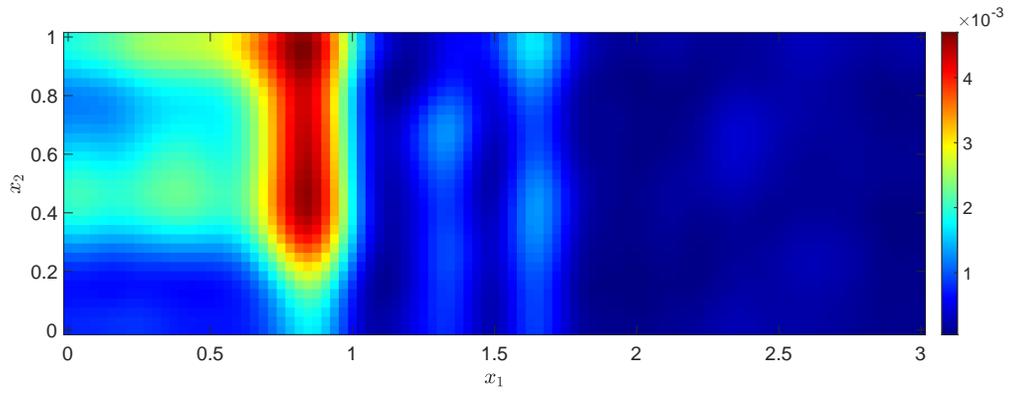} \\
    (a) Variance $\checkV\left[\hata(x,\hatxi)\right]$, DD-MCMC (assembled).\\
	\vspace{0.05mm} \\
    \includegraphics[width=0.80\textwidth]{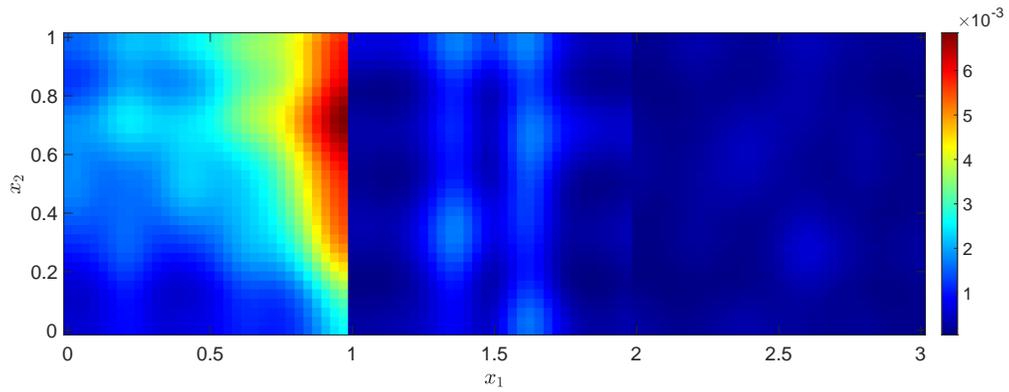}\\
     (b) Variance $\checkV\left[\breve{a}(x,\xi)\right]$, DD-MCMC (stitched).\\
	 \vspace{0.05mm} \\
	 \includegraphics[width=0.80\textwidth]{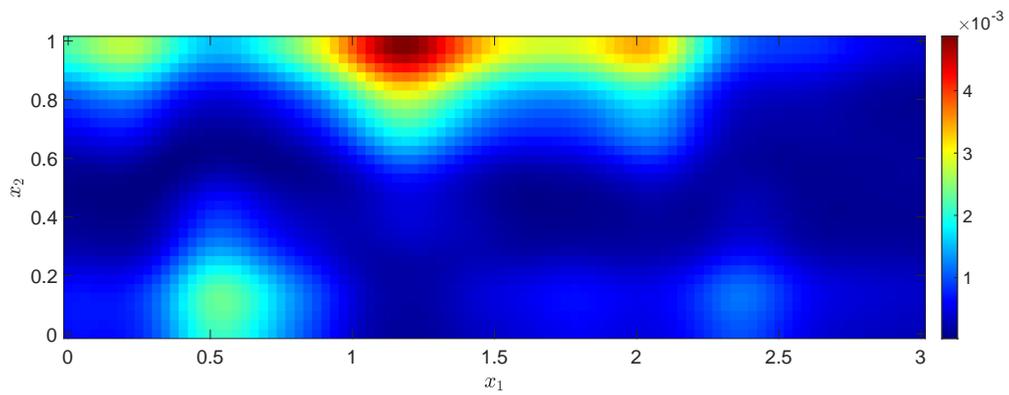}\\
     (c) Variance $\checkV\left[ a(x,\xi) \right]$, G-MCMC.
    \end{tabular}}
        \caption{Estimated variance fields for test problem one ($L=2$).}
    \label{fig:var_compare_L_2}
\end{figure}

\begin{figure}[!htp]
    \centerline{
    \begin{tabular}{c}
    \includegraphics[width=0.80\textwidth]{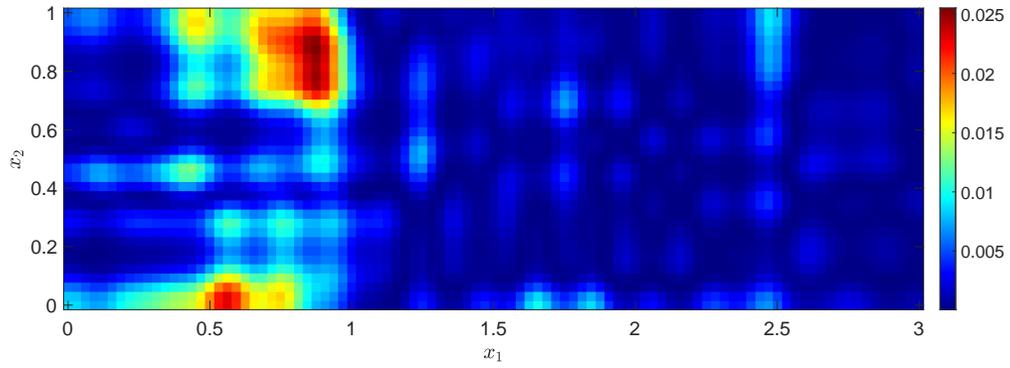} \\
    (a) Variance $\checkV\left[\hata(x,\hatxi)\right]$, DD-MCMC (assembled).\\
	\vspace{0.05mm} \\
    \includegraphics[width=0.80\textwidth]{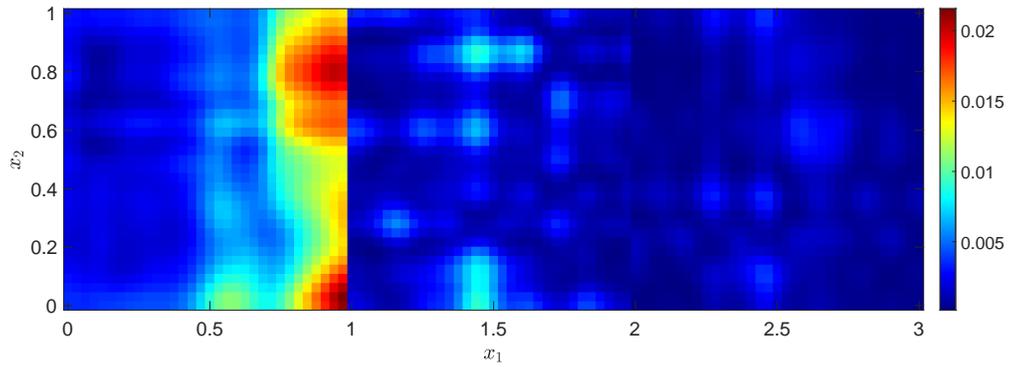}\\
     (b) Variance $\checkV\left[\breve{a}(x,\xi)\right]$, DD-MCMC (stitched).\\
	 \vspace{0.05mm} \\
	 \includegraphics[width=0.80\textwidth]{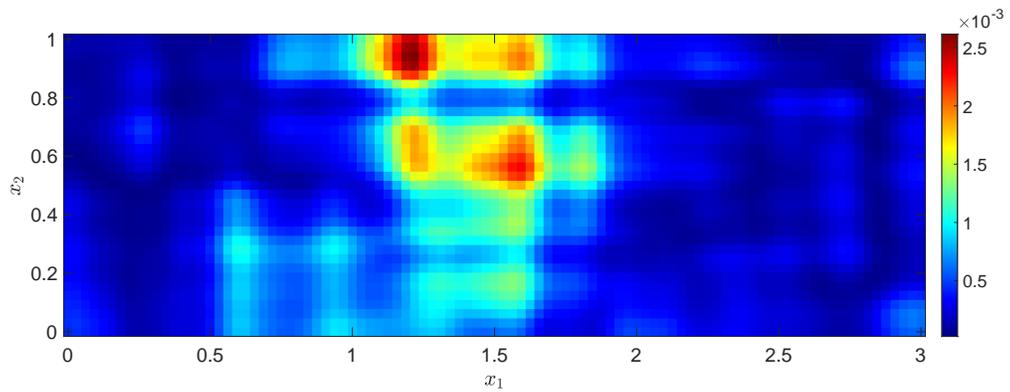}\\
     (c) Variance $\checkV\left[ a(x,\xi) \right]$, G-MCMC.
    \end{tabular}}
    \caption{Estimated variance fields for test problem two ($L=1$).}
    \label{fig:var_compare_L_1}
\end{figure}

\begin{figure}[!htp]
    \centerline{
    \begin{tabular}{c}
    \includegraphics[width=0.80\textwidth]{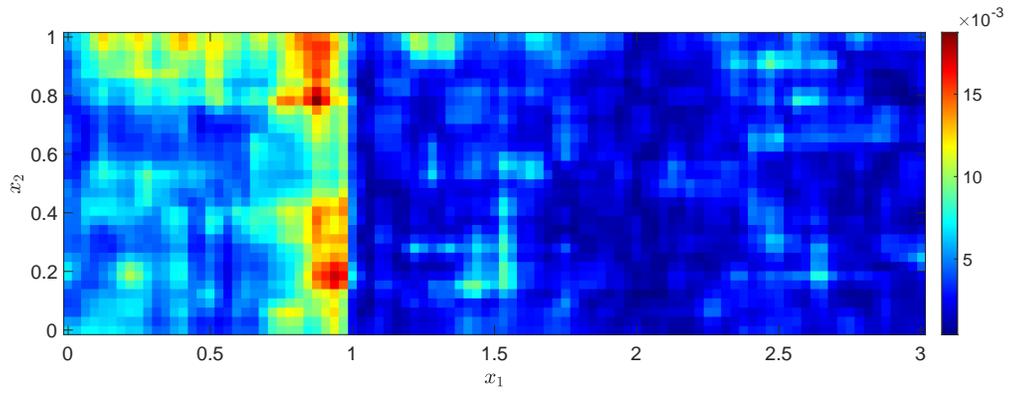} \\
    (a) Variance $\checkV\left[\hata(x,\hatxi)\right]$, DD-MCMC (assembled).\\
	\vspace{0.05mm} \\
    \includegraphics[width=0.80\textwidth]{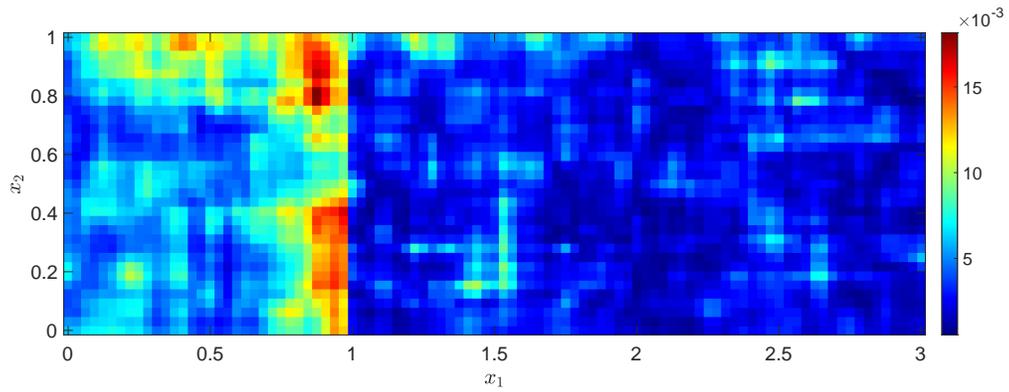}\\
     (b) Variance $\checkV\left[\breve{a}(x,\xi)\right]$, DD-MCMC (stitched).\\
	 \vspace{0.05mm} \\
	 \includegraphics[width=0.80\textwidth]{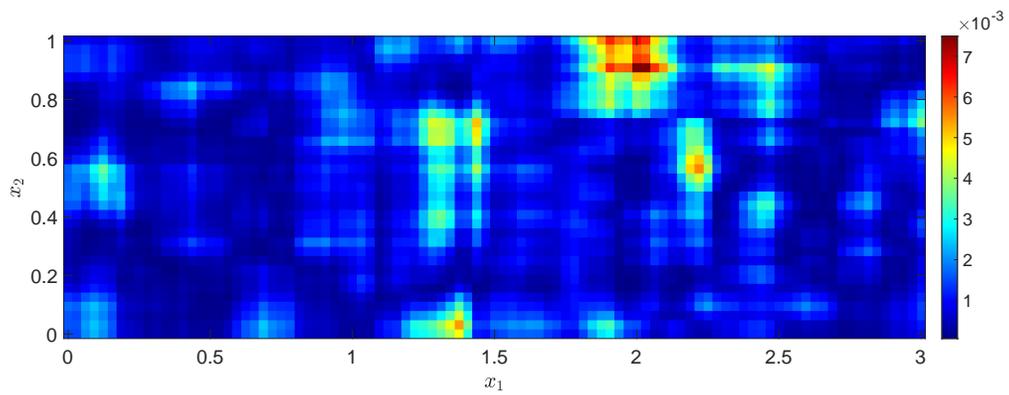}\\
     (c) Variance $\checkV\left[ a(x,\xi) \right]$, G-MCMC.
    \end{tabular}}
    \caption{Estimated variance fields for test problem three ($L=0.5$).}
    \label{fig:var_compare_L_05}
\end{figure}

% posterior mean
For the three test problems, 
Figure \ref{fig:mean_compare_L_2}, Figure \ref{fig:mean_compare_L_1} and Figure \ref{fig:mean_compare_L_05} show the mean fields estimated using the samples obtained from DD-MCMC and G-MCMC. From Figure \ref{fig:mean_compare_L_2}(a), Figure \ref{fig:mean_compare_L_1}(a) and Figure \ref{fig:mean_compare_L_05}(a), it is clear that the estimated field using the DD-MCMC  outputs (the assembled fields $\{\hata(x,\hatxi^{s})\}_{s=1}^{\NC}$) are very similar to the truth permeability fields shown in Figure \ref{fig:setup_2}, Figure \ref{fig:setup_1} and Figure \ref{fig:setup_05}. For test problem one where its truth permeability field is relatively smooth, although the sample mean of the assembled fields gives an accurate estimate (see Figure \ref{fig:mean_compare_L_2}(a)), the sample mean of the stitched fields (i.e.\ $\{\brea(x,\xi^{s}),s=1,\ldots,N\}$) gives misleading information on the interfaces (see Figure \ref{fig:mean_compare_L_2}(b)). This is consistent with our analysis in Section \ref{sec:recon}, and confirms that our reconstruction procedure (line 21 of Algorithm \ref{alg:parallel_MCMC}) is necessary. The results associated with the stitched fields for test problem two and three are shown Figure \ref{fig:mean_compare_L_1}(b) and Figure \ref{fig:mean_compare_L_05}(b). It can be seen that as the truth permeability field of the underlying problem becomes less smooth, the effect of interface becomes less significant. However, it is still clear that the assembled fields give more accuracy mean estimates than the stitched fields. The results of G-MCMC are shown in Figure \ref{fig:mean_compare_L_2}(c), Figure \ref{fig:mean_compare_L_1}(c) and Figure \ref{fig:mean_compare_L_05}(c), where it is clear that for a comparable computational cost, the estimated mean field from G-MCMC is less accurate than that of DD-MCMC (with the assembled fields). In addition, the estimated variance fields  using the samples obtained from DD-MCMC and G-MCMC are shown in 
Figure \ref{fig:var_compare_L_2}, Figure \ref{fig:var_compare_L_1} and Figure \ref{fig:var_compare_L_05}, where it can be seen that the variances are small. 

%the comparisons of the assembled posterior mean field, the stitched posterior mean field and the global inversion. Compared with the stitched posterior mean field, we can see that the assembled field can avoid rough interface conditions. Using the approximately equivalent computational cost, our proposed DD-MCMC algorithm can almost exactly capture the groundtruth while the global inversion fails to depict the actual peamibility field.

%The posterior uncertainties of DD-MCMC and G-MCMC are compared for all three test problems. Figure \ref{fig:var_compare_L_2}, Figure \ref{fig:var_compare_L_1} and Figure \ref{fig:var_compare_L_05} show the comparisons of the posterior variance fileds.  In the $L=1$ and $L=0.5$ cases, the variance fields of the global inversion are typically smaller. However, the relative errors of the global casses are larger, which indicates that the Markov chain generated by the G-MCMC algorithm is probably trapped in a local mode and does not converge to a global optimum.

To assess the accuracy of the estimated posterior mean permeability, we introduce the following quantities of errors
%To access the accuracy of the posterior mean fields, we define the posterior mean relative errors of the global domain and the local domains as 
\begin{align}
& \epsilon:=\| \checkE[a(x,\xi)] -\atruth \|_2/\|\atruth\|_2 ,
\label{eq:re_global}\\
& \breve{\epsilon}:= \| \checkE[\breve{a}(x,\xi)] -\atruth \|_2/\|\atruth\|_2,\label{eq:re_stitch}\\
& \widehat{\epsilon}:=\| \checkE[\hata(x,\hatxi)] -\atruth \|_2/\|\atruth\|_2,\label{eq:re_assemble}
%& \epsilon^{(i)}:=\| \checkE[a^{(i)}(x,\xi^{(i)})] -a^{(i)}_{\text{truth}} \|_2/\|a^{(i)}_{\text{truth}}\|_2\,,\quad i=1,\ldots,M\,, \label{eq:re_local}
\end{align}
where $\atruth$ is the truth permeability, $\checkE[a(x,\xi)]$ is the mean estimate using the samples obtained from G-MCMC, $\checkE[\breve{a}(x,\xi)]$ is the mean estimate using the stitched fields obtained in  DD-MCMC, and $\checkE[\hata(x,\hatxi)]$ is the mean estimate using the assembled fields obtained in  DD-MCMC. Table \ref{tab:rel_err} shows these errors in the mean estimates for the three test problems. It is clear that the errors of our DD-MCMC (for both $\checkE[\hata(x,\hatxi)]$ and $\checkE[\breve{a}(x,\xi)]$) are smaller than the errors of G-MCMC for the three test problems. In addition, the error for the stitched field ($\checkE[\breve{a}(x,\xi)]$) is slightly larger than that for the assembled field ($\checkE[\hata(x,\hatxi)]$), which is consistent with our analysis discussed in Theorem \ref{them:projection}.

\begin{table}[!htp]
	\caption{Errors in mean estimates for the three test problems.}
	\centering
	\begin{tabular}{c|ccc}
		\hline
		$L$&$2$  &$1$  & $0.5$ \\
		\hline
%		$\epsilon^{(1)}$ & $1.084\times 10^{-1}$ & $1.050\times 10^{-1}$ &  $1.355\times 10^{-1}$\\
%		$\epsilon^{(2)}$ &  $2.589\times 10^{-2}$& $8.636\times 10^{-2}$ & $1.035\times 10^{-1}$\\
%		$\epsilon^{(3)}$& $2.234\times 10^{-2}$ & $6.710\times 10^{-2}$ & $9.350\times 10^{-2}$ \\
		$\widehat{\epsilon}$& $5.241\times 10^{-2}$ & $ 7.928\times 10^{-2}$ & $1.083\times 10^{-1}$ \\
		$\breve{\epsilon}$& $5.261\times 10^{-2}$& $8.571\times 10^{-2}$ &  $1.088\times 10^{-1}$\\
		$\epsilon$& $1.340\times 10^{-1}$ & $2.141\times 10^{-1}$&  $1.805\times 10^{-1}$\\
		\hline
	\end{tabular}
	\label{tab:rel_err}
\end{table}

%The convergence properties of DD-MCMC and G-MCMC for three correlation lengths are also studied. Figure \ref{fig:rmse} shows the trends of $\RMSE$ (see \eqref{eq:rmse_global}--\eqref{eq:rmse_local}) against the computational costs in global and local cases. Here the relative errors are computed with all the intermediate posterior samples and the burn-in step is not taken. In all cases, we can see that the $\RMSE$ in the global domain continues to decrease, indicating the convergence of the MCMC algorithm, the global posterior relative error stays at a relatively high place. On the other hand, compared to G-MCMC, we can see that our proposed DD-MCMC algorithm converges faster and can achieve lower recovery accuracy given the same computational cost. In the perspetive of convergence rates and misfit precision, we can see our method outperforms the global inversion. 

%\begin{figure}[!htp]
%    \centerline{
%    \begin{tabular}{ccc}
%    \includegraphics[width=0.32\textwidth]{Figures/Test_Problem1/RMSE_obs.eps} & \includegraphics[width=0.32\textwidth]{Figures/Test_Problem2/RMSE_obs.eps} & \includegraphics[width=0.32\textwidth]{Figures/Test_Problem3/RMSE_obs.eps} \\
%    (a) RMSE$(L=2)$ & (b) RMSE$(L=1)$ & (c) RMSE$(L=0.5)$
%    \end{tabular}}
%    \caption{The convergence results of three test problems. The RMSE of observable pressure values versus computational costs in subdomains and the global domain.}
%    \label{fig:rmse}
%\end{figure}
\section{Conclusion}
\label{sec:conclusion}
The divide and conquer principle is one of the fundamental concepts to solve high-dimensional Bayesian inverse problems involving forward models governed by PDEs.
With a focus on Karhunen-Lo\`{e}ve (KL) expansion based priors, this paper proposes a domain-decomposed Markov chain Monte Carlo (DD-MCMC) algorithm. 
In DD-MCMC, difficulties caused by global prior fields with short correlation lengths are curbed through decomposing global spatial domains into small local subdomains, where correlation lengths become relatively large. On each subdomain, local KL expansion is conducted to result in relatively low-dimensional parameterization, and effective Gaussian process (GP)  interface models are built with an active learning procedure. The global high-dimensional Bayesian inverse problem is then decomposed into a series of local low dimensional-problems, where the corresponding local forward PDE models are also significantly cheaper than the global forward PDE model.  After that, MCMC is applied for local problems to  generate posterior samples of the local input fields. With posterior samples of the local problems, a novel projection procedure is developed to reconstruct samples of the global input field, which are referred to as the assembled fields. Numerical results demonstrate the overall efficiency of the proposed DD-MCMC algorithm. As sufficient data are crucial for efficient Bayesian inversion, although domain decomposition can significantly reduce the computational costs for the inference procedure, small subdomains often only contain limited observation data, which limits the sizes of local subdomains. To further decompose the subdomains, a possible solution is to introduce new effective interface conditions with the whole observed data. Implementing such strategies will be the focus of our future work.

\bigskip
\textbf{Acknowledgments:}
This work is supported by the National Natural Science Foundation of China (No. 12071291), the Science and
Technology Commission of Shanghai Municipality (No. 20JC1414300) and the Natural Science Foundation of Shanghai (No. 20ZR1436200).
%\section*{References}

\bibliography{xu}

\begin{thebibliography}{10}
\expandafter\ifx\csname url\endcsname\relax
  \def\url#1{\texttt{#1}}\fi
\expandafter\ifx\csname urlprefix\endcsname\relax\def\urlprefix{URL }\fi
\expandafter\ifx\csname href\endcsname\relax
  \def\href#1#2{#2} \def\path#1{#1}\fi

\bibitem{springer2021efficient}
S.~Springer, H.~Harrio, J.~Susiluoto, A.~Bibov, A.~Davis, Y.~Marzouk, Efficient
  {Bayesian} inference for large chaotic dynamical systems, Geoscientific Model
  Development 14~(7) (2021) 4319--4333.

\bibitem{solonen2012efficient}
A.~Solonen, P.~Ollinaho, M.~Laine, H.~Haario, J.~Tamminen, H.~J{\"a}rvinen,
  Efficient {MCMC} for climate model parameter estimation: parallel adaptive
  chains and early rejection, Bayesian Analysis 7~(3) (2012) 715--736.

\bibitem{martin2012stochastic}
J.~Martin, L.~C. Wilcox, C.~Burstedde, O.~Ghattas, A stochastic {Newton MCMC}
  method for large-scale statistical inverse problems with application to
  seismic inversion, SIAM Journal on Scientific Computing 34~(3) (2012)
  A1460--A1487.

\bibitem{haario2004markov}
H.~Haario, M.~Laine, M.~Lehtinen, E.~Saksman, J.~Tamminen, {Markov} chain
  {Monte Carlo} methods for high dimensional inversion in remote sensing,
  Journal of the Royal Statistical Society: series B (statistical methodology)
  66~(3) (2004) 591--607.

\bibitem{stuart2010inverse}
A.~M. Stuart, Inverse problems: a {Bayesian} perspective, Acta numerica 19
  (2010) 451--559.

\bibitem{wang04bayesian}
J.~Wang, N.~Zabaras, A {Bayesian} inference approach to the inverse heat
  conduction problem, International Journal of Heat and Mass Transfer 47~(17)
  (2004) 3927--3941.

\bibitem{efendiev2006preconditioning}
Y.~Efendiev, T.~Hou, W.~Luo, Preconditioning {Markov} chain {Monte Carlo}
  simulations using coarse-scale models, SIAM Journal on Scientific Computing
  28~(2) (2006) 776--803.

\bibitem{marzouk2007stochastic}
Y.~M. Marzouk, H.~N. Najm, L.~A. Rahn, Stochastic spectral methods for
  efficient {Bayesian} solution of inverse problems, Journal of Computational
  Physics 224~(2) (2007) 560--586.

\bibitem{lieberman2010parameter}
C.~Lieberman, K.~Willcox, O.~Ghattas, Parameter and state model reduction for
  large-scale statistical inverse problems, SIAM Journal on Scientific
  Computing 32~(5) (2010) 2523--2542.

\bibitem{elmoselhy12bayesian}
T.~A. {El Moselhy}, Y.~M. Marzouk, {Bayesian} inference with optimal maps,
  Journal of Computational Physics 231~(23) (2012) 7815--7850.

\bibitem{tan13computational}
T.~Bui-Thanh, O.~Ghattas, J.~Martin, G.~Stadler, A computational framework for
  infinite-dimensional {Bayesian} inverse problems part {I}: The linearized
  case, with application to global seismic inversion, SIAM Journal on
  Scientific Computing 35~(6) (2013) A2494--A2523.

\bibitem{li2014adaptive}
J.~Li, Y.~M. Marzouk, Adaptive construction of surrogates for the {Bayesian}
  solution of inverse problems, SIAM Journal on Scientific Computing 36~(3)
  (2014) A1163--A1186.

\bibitem{cui16dimension}
T.~Cui, K.~J. Law, Y.~M. Marzouk, Dimension-independent likelihood-informed
  {MCMC}, Journal of Computational Physics 304 (2016) 109--137.

\bibitem{zahm2022certified}
O.~Zahm, T.~Cui, K.~Law, A.~Spantini, Y.~Marzouk, Certified dimension reduction
  in nonlinear {Bayesian} inverse problems, Mathematics of Computation 91~(336)
  (2022) 1789--1835.

\bibitem{robert2013monte}
C.~Robert, G.~Casella, {Monte Carlo} statistical methods, Springer Science \&
  Business Media, 2013.

\bibitem{cui2015data}
T.~Cui, Y.~M. Marzouk, K.~E. Willcox, Data-driven model reduction for the
  {Bayesian} solution of inverse problems, International Journal for Numerical
  Methods in Engineering 102~(5) (2015) 966--990.

\bibitem{chen15sparse}
P.~Chen, C.~Schwab, Sparse-grid, reduced-basis {Bayesian} inversion, Computer
  Methods in Applied Mechanics and Engineering 297 (2015) 84--115.

\bibitem{jiang17multiscale}
L.~Jiang, N.~Ou, Multiscale model reduction method for {Bayesian} inverse
  problems of subsurface flow, Journal of Computational and Applied Mathematics
  319 (2017) 188--209.

\bibitem{liao2019adaptive}
Q.~Liao, J.~Li, An adaptive reduced basis {ANOVA} method for high-dimensional
  {Bayesian} inverse problems, Journal of Computational Physics 396 (2019)
  364--380.

\bibitem{ghanem2003stochastic}
R.~G. Ghanem, P.~D. Spanos, Stochastic finite elements: a spectral approach,
  Courier Corporation, 2003.

\bibitem{li2015note}
J.~Li, A note on the {Karhunen–Lo\`{e}ve} expansions for infinite-dimensional
  {Bayesian} inverse problems, Statistics \& Probability Letters 106 (2015)
  1--4.

\bibitem{ellam2016bayesian}
L.~Ellam, N.~Zabaras, M.~Girolami, A {Bayesian} approach to multiscale inverse
  problems with on-the-fly scale determination, Journal of Computational
  Physics 326 (2016) 115--140.

\bibitem{xia2021bayesian}
Y.~Xia, N.~Zabaras, Bayesian multiscale deep generative model for the solution
  of high-dimensional inverse problems, Journal of Computational Physics 455
  (2022) 111008.

\bibitem{chen2015local}
Y.~Chen, J.~Jakeman, C.~Gittelson, D.~Xiu, Local polynomial chaos expansion for
  linear differential equations with high dimensional random inputs, SIAM
  Journal on Scientific Computing 37~(1) (2015) A79--A102.

\bibitem{liao2015domain}
Q.~Liao, K.~Willcox, A domain decomposition approach for uncertainty analysis,
  SIAM Journal on Scientific Computing 37~(1) (2015) A103--A133.

\bibitem{contreras2018parallela}
A.~A. Contreras, P.~Mycek, O.~P. Le~Ma{\^\i}tre, F.~Rizzi, B.~Debusschere,
  O.~M. Knio, Parallel domain decomposition strategies for stochastic elliptic
  equations. part a: Local {Karhunen--Lo{\`e}ve} representations, SIAM Journal
  on Scientific Computing 40~(4) (2018) C520--C546.

\bibitem{contreras2018parallelb}
A.~A. Contreras, P.~Mycek, O.~P. Le~Ma{\^\i}tre, F.~Rizzi, B.~Debusschere,
  O.~M. Knio, Parallel domain decomposition strategies for stochastic elliptic
  equations part b: Accelerated {Monte} {Carlo} sampling with local {PC}
  expansions, SIAM Journal on Scientific Computing 40~(4) (2018) C547--C580.

\bibitem{khajehpour13domain}
S.~Khajehpour, M.~Hematiyan, L.~Marin, A domain decomposition method for the
  stable analysis of inverse nonlinear transient heat conduction problems,
  International Journal of Heat and Mass Transfer 58~(1) (2013) 125--134.

\bibitem{jagtap2020conservative}
A.~D. Jagtap, E.~Kharazmi, G.~E. Karniadakis, Conservative physics-informed
  neural networks on discrete domains for conservation laws: Applications to
  forward and inverse problems, Computer Methods in Applied Mechanics and
  Engineering 365 (2020) 113028.

\bibitem{jagtap2020extended}
A.~D. Jagtap, G.~E. Karniadakis, Extended physics-informed neural networks
  {(XPINNs)}: A generalized space-time domain decomposition based deep learning
  framework for nonlinear partial differential equations, Communications in
  Computational Physics 28~(5) (2020) 2002--2041.

\bibitem{shukla2021parallel}
K.~Shukla, A.~D. Jagtap, G.~E. Karniadakis, Parallel physics-informed neural
  networks via domain decomposition, Journal of Computational Physics 447
  (2021) 110683.

\bibitem{AINSWORTH19971}
M.~Ainsworth, J.~Oden, A posteriori error estimation in finite element
  analysis, Wiley, 2000.

\bibitem{elman14finite}
H.~Elman, D.~Silvester, A.~Wathen, {Finite Elements} and {Fast Iterative
  Solvers}: with {Applications} in {Incompressible Fluid Dynamics}, Oxford
  University Press (UK), 2014.

\bibitem{metropolis1953equation}
N.~Metropolis, A.~W. Rosenbluth, M.~N. Rosenbluth, A.~H. Teller, E.~Teller,
  Equation of state calculations by fast computing machines, The journal of
  chemical physics 21~(6) (1953) 1087--1092.

\bibitem{hastings1970monte}
W.~K. Hastings, {Monte Carlo} sampling methods using {Markov} chains and their
  applications, Biometrika 57~(1) (1970) 97 -- 109.

\bibitem{le2010spectral}
O.~Le~Ma{\^\i}tre, O.~M. Knio, Spectral methods for uncertainty quantification:
  with applications to computational fluid dynamics, Springer Science \&
  Business Media, 2010.

\bibitem{quarteroni1999domain}
A.~M. Quarteroni, A.~Valli, Domain decomposition methods for partial
  differential equations, Oxford University Press, 1999.

\bibitem{rasmussen2006gaussian}
C.~E. Rasmussen, C.~K. Williams, {Gaussian} process for machine learning, The
  MIT Press, 2006.

\bibitem{rasmussen2010gaussian}
C.~E. Rasmussen, H.~Nickisch, Gaussian processes for machine learning {(GPML)}
  toolbox, The Journal of Machine Learning Research 11 (2010) 3011--3015.

\bibitem{ifiss}
D.~Silvester, H.~Elman, A.~Ramage, {I}ncompressible {F}low and {I}terative
  {S}olver {S}oftware ({IFISS}) version 3.5, {\tt
  http://www.manchester.ac.uk/ifiss/} (September 2016).

\bibitem{roberts2001optimal}
G.~O. Roberts, J.~S. Rosenthal, Optimal scaling for various
  {Metropolis-Hastings} algorithms, Statistical science 16~(4) (2001) 351--367.

\end{thebibliography}

\end{document}